\newcommand{\mytitle}{Games and Ramsey-like cardinals}
\newcommand{\myauthor}{Dan Saattrup Nielsen and Philip Welch}

\documentclass[a4paper,10pt]{article}
	\usepackage[T1]{fontenc}
  \usepackage{amsmath, amssymb, amsfonts} 
  \usepackage{amsthm} 
  \usepackage[latin1]{inputenc} 
  \usepackage{sectsty} 
  \usepackage[light]{antpolt} 
  \usepackage{hyperref} 
  \usepackage{cleveref}	
  \usepackage{graphicx}	
  \usepackage{comment} 
  \usepackage{twoopt} 
  \usepackage{setspace} 
  \usepackage{enumitem}	
  \usepackage{amsbsy} 
  \usepackage{apacite} 
  \usepackage[sectionbib]{natbib} 
	\usepackage{todonotes} 
	\usepackage{wrapfig} 
	
  \usepackage{mathrsfs}	
  \usepackage{stmaryrd}	
  \usepackage{mathabx} 
  \usepackage{proof} 
  \usepackage{tikz} 
    \usetikzlibrary{cd} 
    \usetikzlibrary{matrix,arrows} 
    \usetikzlibrary{decorations.pathreplacing,calc,arrows.meta} 
  \usepackage[all]{xy} 
		\newdir{|>}{-<0pt,5pt>{\blacktriangleup}}
  \usepackage{pigpen} 
  \usepackage{xfrac} 
		
	\usepackage[a4paper, total={6in, 8in}]{geometry}
  \usepackage{wallpaper} 
  \usepackage{fancyhdr}	
  \usepackage{changepage} 

  \usepackage{todonotes} 

\tikzset{
  tree/.style 2 args={
    decorate,
    decoration={
      show path construction,
      lineto code={
        \draw[dotted,-] (\tikzinputsegmentfirst) --($(\tikzinputsegmentfirst)!.5!(\tikzinputsegmentlast)$);
        \draw[-{Latex}] ($(\tikzinputsegmentfirst)!.5!(\tikzinputsegmentlast)$) --(\tikzinputsegmentlast) node [midway,right] {\small{$#1$}};
        \draw[-] (\tikzinputsegmentfirst) --++ (105:0.65cm);
        \draw[-] (\tikzinputsegmentfirst) --++ (75:0.65cm) node [midway, right] {\small{$#2$}};
      }
    }
  }
}

\tikzset{
  treeplain/.style 2 args={
    decorate,
    decoration={
      show path construction,
      lineto code={
        \draw[dotted,-] (\tikzinputsegmentfirst) --($(\tikzinputsegmentfirst)!.5!(\tikzinputsegmentlast)$);
        \draw[-] ($(\tikzinputsegmentfirst)!.5!(\tikzinputsegmentlast)$) --(\tikzinputsegmentlast) node [midway,right] {\small{$#1$}};
        \draw[-] (\tikzinputsegmentfirst) --++ (105:0.65cm);
        \draw[-] (\tikzinputsegmentfirst) --++ (75:0.65cm) node [midway, right] {\small{$#2$}};
      }
    }
  }
}

\setlist{nolistsep}

\newtheoremstyle{scthmstyle} 
	{15pt} 
	{15pt} 
	{\itshape} 
	{} 
	{\bfseries\scshape} 
	{.} 
	{.5em} 
	{} 

\newtheoremstyle{scdefstyle} 
	{15pt} 
	{15pt} 
	{\normalfont} 
	{} 
	{\bfseries\scshape} 
	{.} 
	{.5em} 
	{} 

\newtheoremstyle{scremstyle} 
	{15pt} 
	{15pt} 
	{\normalfont} 
	{} 
	{\itshape} 
	{.} 
	{.5em} 
	{} 

\newtheoremstyle{scclaistyle} 
	{15pt} 
	{15pt} 
	{\normalfont} 
	{0.5cm} 
	{\itshape} 
	{.} 
	{.5em} 
	{} 

\theoremstyle{scthmstyle}
\newtheorem{theorem}[subsection]{Theorem}
\newtheorem{proposition}[subsection]{Proposition}
\newtheorem{lemma}[subsection]{Lemma}
\newtheorem{corollary}[subsection]{Corollary}
\theoremstyle{scdefstyle}
\newtheorem{definition}[subsection]{Definition}
\newtheorem{example}[subsection]{Example}
\newtheorem{question}[subsection]{Question}
\theoremstyle{scremstyle}
\newtheorem{remark}[subsection]{Remark}
\theoremstyle{scclaistyle}
\newtheorem{claim}[subsubsection]{Claim}

  \newcommand{\eq}[1]{\begin{align*} #1 \end{align*}}

  \newcommand{\cd}[1]{\eq{\xymatrix@L=6pt{#1}}}
  
  \renewcommand{\abstract}[1]{\begin{quote}{\footnotesize\textsc{Abstract.} #1}\\\end{quote}}

  \newcommandtwoopt{\theo}[3][][]{
    \begin{theorem}[#1]\label[theorem]{#2}
      #3
    \end{theorem}}
  \newcommandtwoopt{\prop}[3][][]{
    \begin{proposition}[#1]\label[proposition]{#2}
      #3
    \end{proposition}}
  \newcommandtwoopt{\lemm}[3][][]{
    \begin{lemma}[#1]\label[lemma]{#2}
      #3
    \end{lemma}}
  \newcommandtwoopt{\coro}[3][][]{
    \begin{corollary}[#1]\label[corollary]{#2}
      #3
    \end{corollary}}
  \newcommandtwoopt{\defi}[3][][]{
    \begin{definition}[#1]\label[definition]{#2}
      #3$\hfill\dashv$
    \end{definition}}
  \newcommandtwoopt{\exam}[3][][]{
    \begin{example}[#1]\label[example]{#2}
      #3
    \end{example}}
  \newcommandtwoopt{\ques}[3][][]{
    \begin{question}[#1]\label[question]{#2}
      #3
    \end{question}}
  \newcommandtwoopt{\rema}[3][][]{
    \begin{remark}[#1]\label[remark]{#2}
      #3
    \end{remark}}

  \newcommandtwoopt{\qtheo}[3][][]{
    \begin{theorem}[#1]\label[theorem]{#2}
      #3$\hfill\dashv$
    \end{theorem}}
  \newcommandtwoopt{\qprop}[3][][]{
    \begin{proposition}[#1]\label[proposition]{#2}
      #3$\hfill\dashv$
    \end{proposition}}
  \newcommandtwoopt{\qlemm}[3][][]{
    \begin{lemma}[#1]\label[lemma]{#2}
      #3$\hfill\dashv$
    \end{lemma}}
  \newcommandtwoopt{\qcoro}[3][][]{
    \begin{corollary}[#1]\label[corollary]{#2}
      #3$\hfill\dashv$
    \end{corollary}}

  \newcommandtwoopt{\defin}[3][][]{
    \begin{definition}[#1]\label[definition]{#2}
      #3
    \end{definition}}

  \renewcommand{\proof}[1]{\textsc{Proof.} #1$\qed$\\}
  \newcommand{\clai}[2][]{\begin{claim}[#1]#2\end{claim}}
  \newcommand{\cproof}[1]{
    \begin{adjustwidth}{0.5cm}{0pt}
      \textsc{Proof of claim.} #1$\hfill\dashv$\\
    \end{adjustwidth}}
  \renewcommand{\qed}{\hfill\blacksquare}

  \DeclareMathOperator{\M}{\mathcal M}
  \DeclareMathOperator{\N}{\mathcal N}

  \DeclareMathOperator{\T}{\mathcal T}

  \DeclareMathOperator{\G}{\mathcal G}
  \DeclareMathOperator{\h}{\mathcal H}

  \DeclareMathOperator{\on}{On}

  \DeclareMathOperator{\ran}{ran}

  \DeclareMathOperator{\id}{id}

  \DeclareMathOperator{\ult}{Ult}

  \DeclareMathOperator{\cof}{cof}

  \DeclareMathOperator{\crit}{crit}

  \DeclareMathOperator{\contr}{\lightning}
  \DeclareMathOperator{\proves}{\vdash}

  \DeclareMathOperator{\forces}{\Vdash}

  \DeclareMathOperator{\restr}{\upharpoonright}

  \DeclareMathOperator{\col}{Col}

  \DeclareMathOperator{\hull}{Hull}
  \DeclareMathOperator{\chull}{cHull}

  \renewcommand{\P}{\mathcal P}

  \renewcommand{\subset}{\subseteq}
  \renewcommand{\supset}{\supseteq}

	\renewcommand{\l}{|}

	\newcommand{\p}{\mathscr P}
	\newcommand{\pistol}{\mathparagraph}

  \newcommand{\game}[8]{\eq{\begin{array}{ccccccccc} \text{I} & #1 && #3 && #5 && #7\\ \text{II} && #2 && #4 && #6 && #8 \end{array}}}

  \newcommand{\bra}[1]{\langle #1\rangle}

  \newcommand{\abs}[1]{\left|#1\right|}

  \newcommand{\pnormal}{\mathrel{\ooalign{$\lneq$\cr\raise.22ex\hbox{$\lhd$}\cr}}}
  \newcommand{\pideal}{\mathrel{\ooalign{$\lneq$\cr\raise.22ex\hbox{$\lhd$}\cr}}}

  \newcommand{\zfc}{\textsf{ZFC}}

  \newcommand{\godel}[1]{\ulcorner #1 \urcorner}

  \newcommand{\po}{\ar@{}[dr]|{\text{\pigpenfont R}}}
  \newcommand{\pb}{\ar@{}[dr]|{\text{\pigpenfont J}}}
  	
  \newcommand{\pinit}{\lhd}	

\begin{document}
	
	\title{\mytitle}
	\author{{\small\textsc{\myauthor}}}
	\date{\today}
	\maketitle

\abstract{
	We generalise the $\alpha$-Ramsey cardinals introduced in \cite{HolySchlicht} for cardinals $\alpha$ to arbitrary ordinals $\alpha$, and answer several questions posed in that paper. In particular, we show that $\alpha$-Ramseys are downwards absolute to the core model $K$ for all $\alpha$ of uncountable cofinality, that strategic $\omega$-Ramsey cardinals are equiconsistent with remarkable cardinals and that strategic $\alpha$-Ramsey cardinals are equiconsistent with measurable cardinals for all $\alpha>\omega$. We also show that the $n$-Ramseys satisfy indescribability properties and use them to provide a game-theoretic characterisation of completely ineffable cardinals, as well as establishing further connections between the $\alpha$-Ramsey cardinals and the Ramsey-like cardinals introduced in \cite{Ramsey1}, \cite{Feng} and \cite{SharpeWelch}.
}

\section{Introduction}

Most of the large cardinals above measurable cardinals can be characterised as the critical points of elementary embeddings $j:V\to\M$, where the strength of the large cardinal notion in question is increased by requiring more closure of the target model $\M$ and more properties of the embedding $j$. In analogy, Ramsey-like cardinals were introduced in \cite{Ramsey1} and \cite{Ramsey2} to be a natural weakening of this concept, being roughly cardinals $\kappa$ that can be characterised as critical points of elementary embeddings $j:\M\to\N$ between $\kappa$-sized $\zfc^-$-models $\M$ and $\N$. Here we then increase our consistency strength by requiring more closure of the domain model $\M$ and more properties of the embedding $j$.

\qquad Implicit work in \cite{Mitchell} and \cite{Kapplications} shows that Ramsey cardinals are precisely of this type, in which the derived measure from $j$ is both weakly amenable and countably complete.\footnote{For a proof of this result see Theorem 1.3 of \cite{Ramsey1}.} The question is then how many of the well-known large cardinals can be characterised in this fashion? \cite{Ramsey1} introduced various Ramsey-like cardinals, whose definitions we will recall in the next section, and recently \cite{HolySchlicht} have introduced a new family of cardinals, called (strategic) $\alpha$-Ramsey cardinals, which have the added feature of having a game-theoretic definition.

\qquad In \cite{HolySchlicht} the (strategic) $\alpha$-Ramseys were considered for $\alpha$ being an infinite cardinal, and in this paper we will expand this definition to any ordinal $\alpha$. Section 3 will cover the finite case which allows us to characterise ineffable-type cardinals and show indescribability properties of these cardinals --- these arguments are based on arguments in \cite{Abramson}.

\qquad Section 4 contains the countable case in which we establish that strategic $\omega$-Ramseys are equiconsistent with Schindler's \textit{remarkable cardinals}, and use this to show that strategic $\omega$-Ramseys are of strictly stronger consistency strength than the $\omega$-Ramseys. We will also consider a hierarchy between $\omega$-Ramsey cardinals and Ramsey cardinals called $(\omega,\alpha)$-Ramsey cardinals, which we will show interleaves with the $\alpha$-iterable cardinals introduced in \cite{Ramsey1}, and lastly show that $(\omega{+}1)$-Ramseys are Ramsey limits of Ramseys and that strategic $(\omega{+}1)$-Ramseys are equiconsistent with a measurable cardinal.

\qquad In section 5 we investigate how the strongly Ramsey and super Ramsey cardinals introduced in \cite{Ramsey1} relate to the $\alpha$-Ramsey cardinals and show that these latter cardinals are downwards absolute to the core model $K$. The last part of this section is dedicated to showing a tight correspondence between strategic $\alpha$-Ramsey cardinals and the $\alpha$-very Ramsey cardinals introduced in \cite{SharpeWelch}, leading to the result that strategic $\omega_1$-Ramsey cardinals are measurable in the core model $K$ below a Woodin cardinal. Section 6 contains an overview of open problems concerning these Ramsey-like cardinals.

\qquad The last section includes two diagrams, showing the relations between all the Ramsey-like cardinals considered in this paper, both in terms of consistency strength and direct implication. A solid line means that the (consistency or direct) implication is ``strict'', in the sense that no proof exists for the implication in the opposite direction, and a dashed line means that we do not know whether the implication is strict or not.

\section{Setting the scene}
\label{sect.settingthescene}

In this section we will recall a handful of definitions concerning Ramsey-like cardinals, as well as define the $\alpha$-Ramsey cardinals for arbitrary ordinals $\alpha$. We start out with the models and measures that we are going to consider.

\defi{
	For a cardinal $\kappa$, a \textbf{weak $\kappa$-model} is a set $\M$ of size $\kappa$ satisfying that $\kappa+1\subset\M$ and $(\M,\in)\models\zfc^-$. If furthermore $\M^{<\kappa}\subset \M$, $\M$ is a \textbf{$\kappa$-model}.\footnote{Note that our (weak) $\kappa$-models do not have to be transitive, in contrast to the models considered in \cite{Ramsey1} and \cite{Ramsey2}. Not requiring the models to be transitive was introduced in \cite{HolySchlicht}.}
}

Recall that $\mu$ is an \textbf{$\M$-measure} if $(\M,\in,\mu)\models\godel{\mu\text{ is a $\kappa$-complete ultrafilter on $\kappa$}}$.

\defi{
	Let $\M$ be a weak $\kappa$-model and $\mu$ an $\M$-measure. Then $\mu$ is
\begin{itemize}
	\item \textbf{weakly amenable} if $x\cap\mu\in\M$ for every $x\in\M$ with $\M$-cardinality $\kappa$;
	\item \textbf{countably complete} if $\bigcap\vec X\neq\emptyset$ for every $\omega$-sequence $\vec X\in{^\omega\mu}$;
	\item \textbf{$\M$-normal} if $(\M,\in,\mu)\models\forall\vec X\in{^\kappa}\mu:\triangle\vec X\in\mu$;
	\item \textbf{genuine} if $|\triangle\vec X|=\kappa$ for every $\kappa$-sequence $\vec X\in{^\kappa\mu}$;
	\item \textbf{normal} if $\triangle\vec X$ is stationary in $\kappa$ for every $\kappa$-sequence $\vec X\in{^\kappa\mu}$;
	\item \textbf{$0$-good}, or simply \textbf{good}, if it has a well-founded ultrapower;
	\item \textbf{$\alpha$-good} for $\alpha>0$ if it is weakly amenable and has $\alpha$-many well-founded iterates.
\end{itemize}
}

Note that a genuine $\M$-measure is $\M$-normal and countably complete, and a countably complete weakly amenable $\M$-measure is $\alpha$-good for all ordinals $\alpha$. We'll use the fact shown in \cite{HolySchlicht} that an $\M$-measure $\mu$ is normal iff $\triangle\vec X$ is stationary for some enumeration $\vec X=\bra{X_\alpha\mid\alpha<\kappa}$ of $\mu$. We are also going to use the following alternative characterisation of weak amenability.

\qprop[Folklore]{
	Let $\M$ be a weak $\kappa$-model, $\mu$ an $\M$-measure and $j:\M\to\N$ the associated ultrapower embedding. Then $\mu$ is weakly amenable if and only if $j$ is \textbf{$\kappa$-powerset preserving}, meaning that $\M\cap\p(\kappa)=\N\cap\p(\kappa)$.
}

The $\alpha$-Ramsey cardinals in \cite{HolySchlicht} are based upon the following game.\footnote{Unless otherwise stated, every game considered will be a game with perfect information between two players I and II. For a formal framework modelling these games, see e.g. \cite{Kanamori}.}

\defi[Holy-Schlicht]{
	For an uncountable cardinal $\kappa=\kappa^{<\kappa}$, a limit ordinal $\gamma\leq\kappa$ and a regular cardinal $\theta>\kappa$ define the game $wfG_\gamma^\theta(\kappa)$ of length $\gamma$ as follows.
	\game{\M_0}{\mu_0}{\M_1}{\mu_1}{\M_2}{\mu_2}{\cdots}{\cdots}

 Here $\M_\alpha\prec H_\theta$ is a $\kappa$-model and $\mu_\alpha$ is a filter for all $\alpha<\gamma$, such that $\mu_\alpha$ is an $\M_\alpha$-measure, the $\M_\alpha$'s and $\mu_\alpha$'s are $\subset$-increasing and $\bra{\M_\xi\mid\xi<\alpha},\bra{\mu_\xi\mid\xi<\alpha}\in\M_\alpha$ for every $\alpha<\gamma$. Letting $\mu:=\bigcup_{\alpha<\gamma}\mu_\alpha$ and $\M:=\bigcup_{\alpha<\gamma}\M_\alpha$, player II wins iff $\mu$ is an $\M$-normal good $\M$-measure.
}

Recall that two games $G_1$ and $G_2$ are \textbf{equivalent} if player I has a winning strategy in $G_1$ iff they have one in $G_2$, and player II has a winning strategy in $G_1$ iff they have one in $G_2$. \cite{HolySchlicht} showed that the games $wfG_\gamma^{\theta_0}(\kappa)$ and $wfG_\gamma^{\theta_1}(\kappa)$ are equivalent for any $\gamma$ with $\cof\gamma\neq\omega$ and any regular $\theta_0,\theta_1>\kappa$. We will be working with a variant of the $wfG_\gamma(\kappa)$ games in which we require less of player I but more of player II. It will turn out that this change of game is innocuous, as Proposition \ref{prop.tildegame} will show that they are equivalent.

\defi[Holy-Schlicht-N.]{
	Let $\kappa=\kappa^{<\kappa}$ be an uncountable cardinal, $\gamma\leq\kappa$ and $\zeta$ ordinals and $\theta>\kappa$ a regular cardinal. Then define the following game $\G_\gamma^\theta(\kappa,\zeta)$ with $(\gamma{+}1)$-many rounds:
	\game{\M_0}{\mu_0}{\M_1}{\mu_1}{\cdots}{\cdots}{\M_\gamma}{\mu_\gamma}

	Here $\M_\alpha\prec H_\theta$ is a weak $\kappa$-model for every $\alpha\leq\gamma$, $\mu_\alpha$ is a normal $\M_\alpha$-measure for $\alpha<\gamma$, $\mu_\gamma$ is an $\M_\gamma$-normal good $\M_\gamma$-measure and the $\M_\alpha$'s and $\mu_\alpha$'s are $\subset$-increasing. For limit ordinals $\alpha\leq\gamma$ we furthermore require that $\M_\alpha=\bigcup_{\xi<\alpha}\M_\xi$, $\mu_\alpha=\bigcup_{\xi<\alpha}\mu_\xi$ and that $\mu_\alpha$ is $\zeta$-good. Player II wins iff they could continue to play throughout all $(\gamma{+}1)$-many rounds.
}

For convenience we will write $\G_\gamma^\theta(\kappa)$ for the game $\G_\gamma^\theta(\kappa,0)$, and $\G_\gamma(\kappa)$ for $\G_\gamma^\theta(\kappa)$ whenever $\cof\gamma\neq\omega$, as again the existence of winning strategies in these games doesn't depend upon a specific $\theta$. Note that we assume that $\kappa=\kappa^{<\kappa}$ is uncountable in the definition of the games that we're considering, so this is a standing assumption throughout the paper, whenever any one of the above two games are considered.

\prop[Holy-Schlicht-N.]{
\label{prop.tildegame}
$\G_\gamma^\theta(\kappa)$, $\G_\gamma^\theta(\kappa,1)$ and $wfG_\gamma^\theta(\kappa)$ are all equivalent for all limit ordinals $\gamma\leq\kappa$, and $\G_\gamma^\theta(\kappa,\zeta)$ is equivalent to $\G_\gamma^\theta(\kappa)$ whenever $\cof\gamma>\omega$ and $\zeta\in\on$.
}
\proof{
	We start by showing the latter statement, so assume that $\cof\gamma>\omega$. Consider now the auxilliary game, call it $\mathcal G$, which is exactly like $\G_\gamma^\theta(\kappa,0)$, but where we also require that $^\omega{\M_\alpha}\subset\M_{\alpha+1}$ and $\bra{\M_\xi\mid\xi\leq\alpha},\bra{\mu_\xi\mid\xi\leq\alpha}\in\M_{\alpha+1}$ for every $\alpha<\gamma$.

	\clai{
		$\mathcal G$ is equivalent to $\G_\gamma^\theta(\kappa)$.
	}

	\cproof{
		If player I has a winning strategy in $\mathcal G$ then they also have one in $\G_\gamma^\theta(\kappa)$, by doing exactly the same. Analogously, if player II has a winning strategy in $\G_\gamma^\theta(\kappa)$ then they also have one in $\mathcal G$. If player I has a winning strategy $\sigma$ in $\G_\gamma^\theta(\kappa)$ then we can construct a winning strategy $\sigma'$ in $\mathcal G$, which is defined as follows. Fix some $\alpha\leq\gamma$ and, writing $\vec\M_\xi:=\bra{\M_\xi\mid\xi\leq\alpha}$ and $\vec\mu_\xi:=\bra{\mu_\xi\mid\xi\leq\alpha}$, we set
		\eq{
			\sigma'(\bra{\M_\xi,\mu_\xi\mid\xi\leq\alpha}):=\hull^{H_\theta}(\sigma(\bra{\M_\xi,\mu_\xi\mid\xi\leq\alpha})\cup{^\omega{\M_\alpha}}\cup\{\vec\M_\xi,\vec\mu_\xi\}),
		}

		i.e. that we're simply throwing in the sequences into our models and making sure that we're still an elementary substructure of $H_\theta$. This new strategy $\sigma'$ is clearly winning. Assuming now that $\tau$ is a winning strategy for player II in $\mathcal G$, we define a winning strategy $\tau'$ for player II in $\G_\gamma^\theta(\kappa)$ by letting $\tau'(\bra{\M_\xi,\mu_\xi\mid\xi\leq\alpha})$ be the result of throwing in the appropriate sequences into the models $\M_\xi$, applying $\tau$ to get a measure, and intersecting that measure with $\M_\alpha$ to get an $\M_\alpha$-measure.
	}

	Now, letting $\M_\gamma$ be the final model of a play of $\mathcal G$, $\cof\gamma>\omega$ implies that any $\omega$-sequence $\vec X\in\M_\gamma$ really is a sequence of elements from some $\M_\xi$ for $\xi<\gamma$, so that $\vec X\in\M_{\xi+1}$ by definition of $\mathcal G$, making $\M_\gamma$ closed under $\omega$-sequences and thus also $\mu_\gamma$ countably complete. Since $\gamma$ is a limit ordinal and the models contain the previous measures and models as elements, the proof of e.g. Theorem 5.6 in \cite{HolySchlicht} shows that $\mu_\gamma$ is also weakly amenable, making it $\zeta$-good for all ordinals $\zeta$.

	\qquad Now we deal with the first statement, so fix a limit ordinal $\gamma$. Firstly $\G_\gamma^\theta(\kappa)$ is equivalent to $\G_\gamma^\theta(\kappa,1)$ as above, since both are equivalent to the auxilliary game $\G$ when $\gamma$ is a limit ordinal. So it remains to show that $\G_\gamma^\theta(\kappa)$ is equivalent to $wfG_\gamma^\theta(\kappa)$. If player I has a winning strategy $\sigma$ in $wfG_\gamma^\theta(\kappa)$ then define a winning strategy $\sigma'$ for player I in $\G_\gamma^\theta(\kappa)$ as
	\eq{
		\sigma'(\bra{\M_\xi,\mu_\xi\mid\xi\leq\alpha}):=\sigma(\bra{\M_0,\mu_0}^\smallfrown\bra{\M_{\xi+1},\mu_{\xi+1}\mid\xi+1\leq\alpha})
	}

	and for limit ordinals $\alpha\leq\gamma$ set $\sigma'(\bra{\M_\xi,\mu_\xi\mid\xi<\alpha}):=\bigcup_{\xi<\alpha}\M_\xi$; i.e. they simply follow the same strategy as in $wfG_\gamma^\theta(\kappa)$ but plugs in unions at limit stages. Likewise, if player II had a winning strategy in $\G_\gamma^\theta(\kappa)$ then they also have a winning strategy in $wfG_\gamma^\theta(\kappa)$, this time just by skipping the limit steps in $\G_\gamma^\theta(\kappa)$.

	\qquad Now assume that player I has a winning strategy $\sigma$ in $\G_\gamma^\theta(\kappa)$ and that player I \textit{doesn't} have a winning strategy in $wfG_\gamma^\theta(\kappa)$. Then define a strategy $\sigma'$ for player I in $wfG_\gamma^\theta(\kappa)$ as follows. Let $s=\bra{\M_\alpha,\mu_\alpha\mid\alpha\leq\eta}$ be a partial play of $wfG_\gamma^\theta(\kappa)$ and let $s'$ be the modified version of $s$ in which we have 'inserted' unions at limit steps, just as in the above paragraph. We can assume that every $\mu_\alpha$ in $s'$ is good and $\M_\alpha$-normal as otherwise player II has already lost and player I can play anything. Now, we want to show that $s'$ is a valid partial play of $\G_\gamma^\theta(\kappa)$. All the models in $s$ are $\kappa$-models, so in particular weak $\kappa$-models.
	
	\clai{
		Every $\mu_\alpha$ in $s'$ is normal.
	}	
	
	\cproof{
		Assume without loss of generality that $\alpha=\eta$. Let player I play any legal response $\M$ to $s$ in $wfG_\gamma^\theta(\kappa)$ (such a response always exists). If player II can't respond then player I has a winning strategy by simply following $s^\cap\bra{\M}$, $\contr$, so player II \textit{does} have a response $\mu$ to $s^\cap\M$. But now the rules of $wfG_\gamma^\theta(\kappa)$ ensures that $\mu_\eta\in\M$, so since
		\eq{
		(\M,\in,\mu)\models\forall\vec X\in{^\kappa}\mu:\godel{\triangle\vec X\text{ is stationary in }\kappa},
		}

		we then also get that $\M\models\godel{\triangle\mu_\eta\text{ is stationary in $\kappa$}}$ since $\mu_\eta\subset\mu$, so elementarity of $\M$ in $H_\theta$ implies that $\triangle\mu_\eta$ really \textit{is} stationary in $\kappa$, making $\mu_\eta$ normal.
	}
	
	This makes $s'$ a valid partial play of $\G_\gamma^\theta(\kappa)$, so we may form the weak $\kappa$-model $\tilde\M_\eta:=\sigma(s')$. Now let $\M_\eta\prec H_\theta$ be a $\kappa$-model with $\tilde\M_\eta\subset\M_\eta$ and $s\in\M_\eta$ and set $\sigma'(s):=\M_\eta$. This defines the strategy $\sigma'$ for player I in $wfG_\gamma^\theta(\kappa)$, which is winning since the winning condition for the two games is the same for $\gamma$ a limit.\footnote{More precisely, that $\sigma$ is winning in $\G_\gamma^\theta(\kappa)$ means that there's a sequence $\bra{f_n:\kappa\to\kappa\mid n<\omega}$ with the $f_n$'s all being elements of the last model $\tilde\M_\gamma$, witnessing the illfoundedness of the ultrapower. But then all these functions will also be elements of the union of the $\M_\alpha$'s, since we ensured that $\M_\alpha\supset\tilde\M_\alpha$ in the construction above, making the ultrapower of $\bigcup_{\alpha<\gamma}\M_\alpha$ by $\bigcup_{\alpha<\gamma}\mu_\alpha$ illfounded as well.}

	\qquad Next, assume that player II has a winning strategy $\tau$ in $wfG_\gamma^\theta(\kappa)$. We recursively define a strategy $\tilde\tau$ for player II in $\G_\gamma^\theta(\kappa)$ as follows. If $\tilde\M_0$ is the first move by player I in $\G_\gamma^\theta(\kappa)$, let $\M_0\prec H_\theta$ be a $\kappa$-model with $\tilde\M_0\subset\M_0$, making $\M_0$ a valid move for player I in $wfG_\gamma^\theta(\kappa)$. Write $\mu_0:=\tau(\bra{\M_0})$ and then set $\tilde\tau(\bra{\tilde\M_0})$ to be $\tilde\mu_0:=\mu_0\cap\tilde\M_0$, which again is normal by the same trick as above, making $\tilde\mu_0$ a legal move for player II in $\G_\gamma^\theta(\kappa)$. Successor stages $\alpha+1$ in the construction are analogous, but we also make sure that $\bra{\M_\xi\mid\xi<\alpha+1},\bra{\mu_\xi\mid\xi<\alpha+1}\in\M_{\alpha+1}$. At limit stages $\tau$ outputs unions, as is required by the rules of $\G_\gamma^\theta(\kappa)$. Since the union of all the $\mu_\alpha$'s is good as $\tau$ is winning, $\tilde\mu_\gamma:=\bigcup_{\alpha<\gamma}\tilde\mu_\alpha$ is good as well, making $\tilde\tau$ winning and we are done.
}

We now arrive at the definitions of the cardinals we will be considering. They were in \cite{HolySchlicht} only defined for $\gamma$ being a cardinal, but given the above result we generalise it to all ordinals $\gamma$.

\defi{
	Let $\kappa$ be a cardinal and $\gamma\leq\kappa$ an ordinal. Then $\kappa$ is \textbf{$\gamma$-Ramsey} if player I does not have a winning strategy in $\G_\gamma^\theta(\kappa)$ for all regular $\theta>\kappa$. We furthermore say that $\kappa$ is \textbf{strategic $\gamma$-Ramsey} if player II \textit{does} have a winning strategy in $\G_\gamma^\theta(\kappa)$ for all regular $\theta>\kappa$. Define \textbf{(strategic) genuine $\gamma$-Ramseys} and \textbf{(strategic) normal $\gamma$-Ramseys} analogously, but where we require the last measure $\mu_\gamma$ to be genuine and normal, respectively.
}

\defi[N.]{
	\label{defi.cohramsey}
	A cardinal $\kappa$ is \textbf{${<}\gamma$-Ramsey} if it is $\alpha$-Ramsey for every $\alpha<\gamma$, \textbf{almost fully Ramsey} if it is ${<}\kappa$-Ramsey and \textbf{fully Ramsey} if it is $\kappa$-Ramsey. Further, say that $\kappa$ is \textbf{coherent ${<}\gamma$-Ramsey} if it's strategic $\alpha$-Ramsey for every $\alpha<\gamma$ and that there exists a choice of winning strategies $\tau_\alpha$ in $\G_\alpha(\kappa)$ for player II satisfying that $\tau_\alpha\subset\tau_\beta$ whenever $\alpha<\beta$. In other words, there is a single strategy $\tau$ for player II in $\G_\gamma(\kappa)$ such that $\tau$ is a winning strategy for player II in $\G_\alpha(\kappa)$ for every $\alpha<\gamma$.\footnote{Note that, with this terminology, ``coherent'' is a stronger notion than ``strategic''. We could've called the cardinals \textit{coherent strategic ${<}\gamma$-Ramseys}, but we opted for brevity instead.}
}

This is not the original definition of (strategic) $\gamma$-Ramsey cardinals however, as this involved elementary embeddings between weak $\kappa$-models -- but as the following theorem of \cite{HolySchlicht} shows, the two definitions coincide whenever $\gamma$ is a regular cardinal.

\qtheo[Holy-Schlicht]{
	\label{theo.Ramseydef}
	For regular cardinals $\lambda$, a cardinal $\kappa$ is $\lambda$-Ramsey iff for arbitrarily large $\theta>\kappa$ and every $A\subset\kappa$ there is a weak $\kappa$-model $\M\prec H_\theta$ with $\M^{<\lambda}\subset\M$ and $A\in\M$ with an $\M$-normal 1-good $\M$-measure $\mu$ on $\kappa$.
}

\section{The finite case}

In this section we are going to consider properties of the $n$-Ramsey cardinals for finite $n$. Note in particular that the $\G_n^\theta(\kappa)$ games are determined, making the ``strategic'' adjective superfluous in this case. We further note that the $\theta$'s are also dispensible in this finite case:

\prop[N.]{
	Let $\kappa<\theta$ be regular cardinals and $n<\omega$. Then player II has a winning strategy in $\G_n^\theta(\kappa)$ iff they have a winning strategy in the game $\G_n(\kappa)$, which is defined as $\G_n^\theta(\kappa)$ except that we don't require that $\M_n\prec H_\theta$.
}
\proof{
	$\Leftarrow$ is clear, so assume that II has a winning strategy $\tau$ in $\G_n^\theta(\kappa)$. Whenever player I plays $\M_k$ in $\G_n(\kappa)$ for $k\leq n$ then define $\M_k^*:=\hull^{H_\theta}(\P)$ where $\P\cong\M_k$ is the transitive collapse of $\M_k$, and play $\M_k^*$ in $\G_n^\theta(\kappa)$. Let $\mu_k$ be the $\tau$-responses to the $\M_k^*$'s and let player II play the $\mu_k$'s in $\G_n(\kappa)$ as well.

	\qquad Assume that this new strategy isn't winning for player II in $\G_n(\kappa)$, so that $\ult(\M_n,\mu_n)$ is illfounded. This is witnessed by some $\omega$-sequence $\vec f:=\bra{f_k\mid k<\omega}$ of $f_k\in{^\kappa o(\M_n)}\cap\M_n$ with $X_k:=\{\alpha<\kappa\mid f_{k+1}(\alpha)<f_k(\alpha)\}\in\mu_n$ for all $k<\omega$. Let $\nu\gg\kappa$, $\h:=\chull^{H_\nu}(\M_n\cup\{\vec f,\M_n,\mu_n\})$ be the transitive collapse of the Skolem hull $\hull^{H_\nu}(\M_n\cup\{\vec f,\M_n,\mu_n\})$, and $\pi:\h\to H_\nu$ be the uncollapse; write $\bar x:=\pi^{-1}(x)$ for all $x\in\ran\pi$.

	\qquad Now $\bar A=A$ for every $A\in\p(\kappa)\cap\M_n$ and thus also $\bar\mu_n=\mu_n$. But now the $\bar f_k$'s witness that $\ult(\bar\M_n,\mu_n)$ is illfounded and thus also that $\ult(\M_n^*,\mu_n)$ is illfounded since $\M_n^*=\hull^{H_\theta}(\bar\M_n)$, contradicting that $\tau$ is winning.
}

For this reason we'll work with the $\G_n(\kappa)$ games throughout this section. Since we don't have to deal with the $\theta$'s anymore we note that $n$-Ramseyness can now be described using a $\Pi^1_{2n+2}$-formula and normal $n$-Ramseyness using a $\Pi^1_{2n+3}$-formula.

\qquad We already have the following characterisations, as proven in \cite{Abramson}.

\theo[Abramson et al.]{
	Let $\kappa=\kappa^{<\kappa}$ be a cardinal. Then
	\begin{enumerate}
		\item $\kappa$ is weakly compact if and only if it is $0$-Ramsey;
		\item $\kappa$ is weakly ineffable if and only if it is genuine $0$-Ramsey;
		\item $\kappa$ is ineffable if and only if it is normal $0$-Ramsey.
	\end{enumerate}
}
\proof{
	This is mostly a matter of changing terminology from \cite{Abramson} to the current game-theoretic one, so we only show $(i)$. Theorem 1.1.3 in \cite{Abramson} shows that $\kappa$ is weakly compact if and only if every $\kappa$-sized collection of subsets of $\kappa$ is measured by a ${<}\kappa$-complete measure, in the sense that every ${<}\kappa$-sequence (in $V$) of measure one sets has non-empty intersection.
	
	\qquad For the $\Rightarrow$ direction we can let player II respond to any $\M_0$ by first getting the ${<}\kappa$-complete $\M_0$-measure $\nu_0$ on $\kappa$ from the above-mentioned result, forming the (well-founded) ultrapower $\pi:\M_0\to\ult(\M_0,\nu)$ and then playing the derived measure of $\pi$, which is $\M_0$-normal and good. For $\Leftarrow$, if $X\subset\p(\kappa)$ has size $\kappa$ then, using that $\kappa=\kappa^{<\kappa}$, we can find a $\kappa$-model $\M_0\prec H_\theta$ with $X\subset\M_0$. Letting player I play $\M_0$ in $\G_0(\kappa)$ we get some $\M_0$-normal good $\M_0$-measure $\mu_0$ on $\kappa$. Since $\M_0$ is closed under ${<}\kappa$-sequences we get that $\mu_0$ is ${<}\kappa$-complete.
}

\subsection*{Indescribability}

In this section we aim to prove that $n$-Ramseys are $\Pi^1_{2n+1}$-indescribable and that normal $n$-Ramseys are $\Pi^1_{2n+2}$-indescribable, which will also establish that the hierarchy of alternating $n$-Ramseys and normal $n$-Ramseys forms a strict hierarchy. Recall the following definition.

\defi{
	A cardinal $\kappa$ is \textbf{$\Pi^1_n$-indescribable} if whenever $\varphi(v)$ is a $\Pi_n$ formula, $X\subset V_\kappa$ and $V_{\kappa+1}\models\varphi[X]$, then there is an $\alpha<\kappa$ such that $V_{\alpha+1}\models\varphi[X\cap V_\alpha]$.
}

Our first indescribability result is then the following, where the $n=0$ case is inspired by the proof of weakly compact cardinals being $\Pi^1_1$-indescribable --- see \cite{Abramson}.

\theo[N.]{
  \label{theo.ind}
	Every $n$-Ramsey $\kappa$ is $\Pi^1_{2n+1}$-indescribable for $n<\omega$.
}
\proof{
  Let $\kappa$ be $n$-Ramsey and assume that it is not $\Pi^1_{2n+1}$-indescribable, witnessed by a $\Pi_{2n+1}$-formula $\varphi(v)$ and a subset $X\subset V_\kappa$, meaning that $V_{\kappa+1}\models\varphi[X]$ and, for every $\alpha<\kappa$, $V_{\alpha+1}\models\lnot\varphi[X\cap V_\alpha]$. We will deal with the $(2n+1)$-many quantifiers occuring in $\varphi$ in $(n+1)$-many steps. We will here describe the first two steps with the remaining steps following the same pattern.

\qquad \framebox{\textbf{First step.}} Write $\varphi(v)\equiv\forall v_1\psi(v,v_1)$ for a $\Sigma_{2n}$-formula $\psi(v,v_1)$. As we are assuming that $V_{\alpha+1}\models\lnot\varphi[X\cap V_\alpha]$ holds for every $\alpha<\kappa$, we can pick witnesses $A^{(0)}_\alpha\subset V_\alpha$ to the outermost existential quantifier in $\lnot\varphi[X\cap V_\alpha]$.

\qquad Let $\M_0$ be a weak $\kappa$-model such that $V_\kappa\subset\M_0$ and $\vec A^{(0)},X\in\M_0$. Fix a good $\M_0$-normal $\M_0$-measure $\mu_0$ on $\kappa$, using the $0$-Ramseyness of $\kappa$. Form $\mathcal A^{(0)}:=[\vec A^{(0)}]_{\mu_0}\in\ult(\M_0,\mu_0)$, where we without loss of generality may assume that the ultrapower is transitive. $\M_0$-normality of $\mu_0$ implies that $\mathcal A^{(0)}\subset V_\kappa$, so that we have that $V_{\kappa+1}\models\psi[X,\mathcal A^{(0)}]$. Now \L o\' s' Lemma, $\M_0$-normality of $\mu_0$ and $V_\kappa\subset\M_0$ also ensures that
\eq{
  \ult(\M_0,\mu_0)\models\godel{V_{\kappa+1}\models\lnot\psi[X,\mathcal A^{(0)}]}.\tag*{$(1)$}
}

This finishes the first step. Note that if $n=0$ then $\lnot\psi$ would be a $\Delta_0$-formula, so that $(1)$ would be absolute to the true $V_{\kappa+1}$, yielding a contradiction. If $n>0$ we cannot yet conclude this however, but that is what we are aiming for in the remaining steps.\\

\qquad \framebox{\textbf{Second step.}} Write $\psi(v,v_1)\equiv\exists v_2\forall v_3\chi(v,v_1,v_2,v_3)$ for a $\Sigma_{2(n-1)}$-formula $\chi(v,v_1,v_2,v_3)$. Since we have established that $V_{\kappa+1}\models\psi[X,\mathcal A^{(0)}]$ we can pick some $B^{(0)}\subset V_\kappa$ such that
\eq{
  V_{\kappa+1}\models\forall v_3\chi[X,\mathcal A^{(0)},B^{(0)},v_3]\tag*{$(2)$}
}

which then also means that, for every $\alpha<\kappa$,
\eq{
  V_{\alpha+1}\models\exists v_3\lnot\chi[X\cap V_\alpha,A^{(0)}_\alpha,B^{(0)}\cap V_\alpha,v_3].\tag*{$(3)$}
}

Fix witnesses $A^{(1)}_\alpha\subset V_\alpha$ to the existential quantifier in $(3)$ and define the sets
\eq{
  S_\alpha^{(0)}:=\{\xi<\kappa\mid A_\xi^{(0)}\cap V_\alpha=\mathcal A^{(0)}\cap V_\alpha\}
}

for every $\alpha<\kappa$ and note that $S_\alpha^{(0)}\in\mu_0$ for every $\alpha<\kappa$, since $V_\kappa\subset\M_0$ ensures that $\mathcal A^{(0)}\cap V_\alpha\in\M_0$ and $\M_0$-normality of $\mu_0$ then implies that $S_\alpha^{(0)}\in\mu_0$ is equivalent to
\eq{
  \ult(\M_0,\mu_0)\models\mathcal A^{(0)}\cap V_\alpha=\mathcal A^{(0)}\cap V_\alpha,
}

which is clearly the case. Now let $\M_1\supset\M_0$ be a weak $\kappa$-model such that $\mathcal A^{(0)},\vec A^{(1)},\vec S^{(0)},B^{(0)}\in\M_1$. Let $\mu_1\supset\mu_0$ be an $\M_1$-normal $\M_1$-measure on $\kappa$, using the $1$-Ramseyness of $\kappa$, so that $\M_1$-normality of $\mu_1$ yields that $\triangle\vec S^{(0)}\in\mu_1$. Observe that $\xi\in\triangle\vec S^{(0)}$ if and only if $A^{(0)}_\xi\cap V_\alpha=\mathcal A^{(0)}\cap V_\alpha$ for every $\alpha<\xi$, so if $\xi$ is a limit ordinal then it holds that $A^{(0)}_\xi=\mathcal A^{(0)}\cap V_\xi$. Now, as before, form $\mathcal A^{(1)}:=[\vec A^{(1)}]_{\mu_1}\in\ult(\M_1,\mu_1)$, so that $(2)$ implies that
\eq{
  V_{\kappa+1}\models\chi[X,\mathcal A^{(0)},B^{(0)},\mathcal A^{(1)}]
}

and the definition of the $A_\alpha^{(1)}$'s along with $(3)$ gives that, for every $\alpha<\kappa$,
\eq{
  V_{\alpha+1}\models\lnot\chi[X\cap V_\alpha,A_\alpha^{(0)},B^{(0)}\cap V_\alpha,A_\alpha^{(1)}].
}

Now this, paired with the above observation regarding $\triangle\vec S^{(0)}$, means that for every $\alpha\in\triangle\vec S^{(0)}\cap\text{Lim}$ we have that
\eq{
  V_{\alpha+1}\models\lnot\chi[X\cap V_\alpha,\mathcal A^{(0)}\cap V_\alpha,B^{(0)}\cap V_\alpha,A_\alpha^{(1)}],
}

so that $\M_1$-normality of $\mu_1$ and \L o\' s' lemma implies that
\eq{
  \ult(\M_1,\mu_1)\models\godel{V_{\kappa+1}\models\lnot\chi[X,\mathcal A^{(0)},B^{(0)},\mathcal A^{(1)}]}.
}

This finishes the second step. Continue in this way for a total of $(n+1)$-many steps, ending with a $\Delta_0$-formula $\phi(v,v_1,\hdots,v_{2n+1})$ such that
\eq{
	V_{\kappa+1}\models\phi[X,\mathcal A^{(0)},B^{(0)},\hdots,\mathcal A^{(n-1)},B^{(n-1)},\mathcal A^{(n)}]\tag*{$(4)$}
}

and that $\ult(\M_n,\mu_n)\models\godel{V_{\kappa+1}\models\lnot\phi[X,\mathcal A^{(0)},B^{(0)},\hdots,\mathcal A^{(n)}]}$. But now absoluteness of $\lnot\phi$ means that $V_{\kappa+1}\models\lnot\phi[X,\mathcal A^{(0)},B^{(0)},\hdots,\mathcal A^{(n)}]$, contradicting $(4)$.
}

Note that this is optimal, as $n$-Ramseyness can be described by a $\Pi^1_{2n+2}$-formula. As a corollary we then immediately get the following.

\qcoro[N.]{
	\label{coro.ind}
	Every ${<}\omega$-Ramsey cardinal is $\Delta^2_0$-indescribable.
}

The second indescribability result concerns the normal $n$-Ramseys, where the $n=0$ case here is inspired by the proof of ineffable cardinals being $\Pi^1_2$-indescribable --- see \cite{Abramson}.

\theo[N.]{
	\label{theo.normind}
	Every normal $n$-Ramsey $\kappa$ is $\Pi^1_{2n+2}$-indescribable for $n<\omega$.
}

Before we commence with the proof, note that we cannot simply do the same thing as we did in the proof of Theorem \ref{theo.ind}, as we would end up with a $\Pi^1_1$ statement in an ultrapower, and as $\Pi^1_1$ statements are not upwards absolute in general we would not be able to get our contradiction.\\

\proof{
	Let $\kappa$ be normal $n$-Ramsey and assume that it is not $\Pi^1_{2n+2}$-indescribable, witnessed by a $\Pi_{2n+2}$-formula $\varphi(v)$ and a subset $X\subset V_\kappa$. Use that $\kappa$ is $n$-Ramsey to perform the same $n+1$ steps as in the proof of Theorem \ref{theo.ind}. This gives us a $\Sigma_1$-formula $\phi(v,v_1,\hdots,v_{2n+1})$ along with sequences $\bra{\mathcal A^{(0)},\cdots,\mathcal A^{(n)}}$, $\bra{B^{(0)},\hdots,B^{(n-1)}}$ and a play $\bra{\M_k,\mu_k\mid k\leq n}$ of $\G_n(\kappa)$ in which player II wins and $\mu_n$ is normal, such that
\eq{
  V_{\kappa+1}\models\phi[X,\mathcal A^{(0)},B^{(0)},\hdots,\mathcal A^{(n-1)},B^{(n-1)},\mathcal A^{(n)}]\tag*{$(1)$}
}

and, for $\mu_n$-many $\alpha<\kappa$,
\eq{
  V_{\alpha+1}\models\lnot\phi[X\cap V_\alpha,\mathcal A^{(0)}\cap V_\alpha,B^{(0)}\cap V_\alpha,\hdots,\mathcal A^{(n-1)}\cap V_\alpha,B^{(n-1)}\cap V_\alpha,A^{(n)}_\alpha].
}

Now form $S^{(n)}_\alpha\in\mu_n$ as in the proof of Theorem \ref{theo.ind}. The main difference now is that we do not know if $\vec S^{(n)}\in\M_n$ (in the proof of Theorem \ref{theo.ind} we only ensured that $\vec S^{(k)}\in\M_{k+1}$ for every $k<n$ and we only defined $\vec S^{(k)}$ for $k<n$), but we can now use normality\footnote{Recall that this is stronger than just requiring it to be $\M_n$-normal --- we don't require $\vec S^{(n)}\in\M_n$.} of $\mu_n$ to ensure that we \textit{do} have that $\triangle\vec S^{(n)}$ is stationary in $\kappa$. This means that we get a stationary set $S\subset\kappa$ such that for every $\alpha\in S$ it holds that
\eq{
  V_{\alpha+1}\models\lnot\phi[X\cap V_\alpha,\mathcal A^{(0)}\cap V_\alpha,B^{(0)}\cap V_\alpha,\hdots, B^{(n-1)}\cap V_\alpha,\mathcal A^{(n)}\cap V_\alpha].\tag*{$(2)$}
}

Now note that since $\kappa$ is inaccessible it is $\Sigma^1_1$-indescribable, meaning that we can reflect $(1)$. Furthermore, Lemma 3.4.3 of \cite{Abramson} shows that the set of reflection points of $\Sigma^1_1$-formulas is in fact club, so intersecting this club with $S$ we get a $\zeta\in S$ satisfying that
\eq{
  V_{\zeta+1}\models\phi[X\cap V_\zeta,\mathcal A^{(0)}\cap V_\zeta,B^{(0)}\cap V_\zeta,\hdots, B^{(n-1)}\cap V_\zeta,\mathcal A^{(n)}\cap V_\zeta],
}

contradicting $(2)$.
}

Note that this is optimal as well, since normal $n$-Ramseyness can be described by a $\Pi^1_{2n+3}$-formula. In particular this then means that every $(n{+}1)$-Ramsey is a normal $n$-Ramsey stationary limit of normal $n$-Ramseys, and every normal $n$-Ramsey is an $n$-Ramsey stationary limit of $n$-Ramseys, making the hierarchy of alternating $n$-Ramseys and normal $n$-Ramseys a strict hierarchy.

\subsection*{Downwards absoluteness to $L$}

The following proof is basically the proof of Theorem 4.1.1 in \cite{Abramson}.

\pagebreak
\theo[N.]{
	Genuine- and normal $n$-Ramseys are downwards absolute to $L$, for every $n<\omega$.
}
\proof{
	Assume first that $n=0$ and that $\kappa$ is a genuine $0$-Ramsey cardinal. Let $\M\in L$ be a weak $\kappa$-model --- we want to find a genuine $\M$-measure inside $L$. By assumption we \textit{can} find such a measure $\mu$ in $V$; we will show that in fact $\mu\in L$. Fix any enumeration $\bra{A_\xi\mid\xi<\kappa}\in L$ of $\p(\kappa)\cap\M$. It then clearly suffices to show that $T\in L$, where $T:=\{\alpha<\kappa\mid A_\xi\in\mu\}$.

  \clai{
		\label{clai.pospart}
		$T\cap\alpha\in L$ for any $\alpha<\kappa$.
	}

	\cproof{
		Let $\vec B$ be the \textbf{$\mu$-positive part} of $\vec A$, meaning that $B_\xi:=A_\xi$ if $A_\xi\in\mu$ and $B_\xi:=\lnot A_\xi$ if $A_\xi\notin\mu$. As $\mu$ is genuine we get that $\triangle\vec B$ has size $\kappa$, so we can pick $\delta\in\triangle\vec B$ with $\delta>\alpha$. Then $T\cap\alpha=\{\xi<\alpha\mid\delta\in A_\xi\}$, which can be constructed within $L$.
	}

	But now Lemma 4.1.2 in \cite{Abramson} shows that there is a $\Pi_1$ formula $\varphi(v)$ such that, given any non-zero ordinal $\zeta$, $V_{\zeta+1}\models\varphi[A]$ if and only if $\zeta$ is a regular cardinal and $A$ is a non-constructible subset of $\zeta$. If we therefore assume that $T\notin L$ then $V_{\kappa+1}\models\varphi[T]$, which by $\Pi^1_1$-indescribability of $\kappa$ means that there exists some $\alpha<\kappa$ such that $V_{\alpha+1}\models\varphi[T\cap V_\alpha]$, i.e. that $T\cap\alpha\notin L$, contradicting the claim. Therefore $\mu\in L$. It is still genuine in $L$ as $(\triangle\mu)^L=\triangle\mu$, and if $\mu$ was normal then that is still true in $L$ as clubs in $L$ are still clubs in $V$. The cases where $\kappa$ is a genuine- or normal $n$-Ramsey cardinal is analogous.
}

Since $(n{+}1)$-Ramseys are normal $n$-Ramseys we then immediately get the following.

\qcoro[N.]{
	Every $(n{+}1)$-Ramsey is normal $n$-Ramsey in $L$, for every $n<\omega$. In particular, ${<}\omega$-Ramseys are downwards absolute to $L$.
}

\subsection*{Complete ineffability}

In this section we provide a characterisation of the \textit{completely ineffable} cardinals in terms of the $\alpha$-Ramseys. To arrive at such a characterisation, we need a slight strengthening of the ${<}\omega$-Ramsey cardinals, namely the \textit{coherent ${<}\omega$-Ramseys} as defined in \ref{defi.cohramsey}. Note that a coherent ${<}\omega$-Ramsey is precisely a cardinal satisfying the $\omega$-filter property, as defined in \cite{HolySchlicht}.

\qquad The following theorem shows that assuming coherency does yield a strictly stronger large cardinal notion. The idea of its proof is very closely related to the proof of Theorem \ref{theo.normind} (the indescribability of normal $n$-Ramseys), but the main difference is that we want everything to occur locally inside our weak $\kappa$-models.

\theo[N.]{
	Every coherent ${<}\omega$-Ramsey is a stationary limit of ${<}\omega$-Ramseys.
}
\proof{
	Let $\kappa$ be coherent ${<}\omega$-Ramsey. Let $\theta\gg\kappa$ be regular and let $\M_0\prec H_\theta$ be a weak $\kappa$-model with $V_\kappa\subset\M_0$. Let then player I play arbitrarily while player II plays according to her coherent winning strategies in $\G_n(\kappa)$, yielding a weak $\kappa$-model $\M\prec H_\theta$ with an $\M$-normal $\M$-measure $\mu:=\bigcup_{n<\omega}\mu_n$ on $\kappa$.
		
	\qquad Assume towards a contradiction that $X:=\{\xi<\kappa\mid\xi\text{ is ${<}\omega$-Ramsey}\}\notin\mu$. Since $X=\bigcap\vec X$ and $\vec X\in\M$, where $X_n:=\{\xi<\kappa\mid\xi\text{ is $n$-Ramsey}\}$, we must have by $\M$-normality of $\mu$ that $\lnot X_k\in\mu$ for some $k<\omega$. Note that $\lnot X_k\in\M_0$ by elementarity, so that $\lnot X_k\in\mu_0$ as well. Perform the $k+1$ steps as in the proof of Theorem \ref{theo.normind} with $\varphi(\xi)$ being $\godel{\xi\text{ is $k$-Ramsey}}$, so that we get a weak $\kappa$-model $\M_{k+1}\prec H_\theta$, an $\M_{k+1}$-normal $\M_{k+1}$-measure $\tilde\mu_{k+1}$ on $\kappa$, a $\Sigma_1$-formula $\varphi(v,v_1,v_2,\hdots,v_{2k+1})$ and sequences $\bra{\mathcal A^{(0)},\hdots,\mathcal A^{(k)}}$ and $\bra{B^{(0)},\hdots,B^{(k-1)}}$ such that
	\eq{
		V_{\kappa+1}\models\varphi[\kappa,\mathcal A^{(0)},B^{(0)},\mathcal A^{(1)},B^{(1)},\hdots,\mathcal A^{(k-1)},B^{(k-1)},\mathcal A^{(k)}]\tag*{$(2)$}
	}

	and there is a $Y\in\tilde\mu_{k+1}$ with $Y\subset\lnot X_k$ such that given any $\xi\in Y$,
	\eq{
		V_{\xi+1}\models\lnot\varphi[\xi,A_\xi^{(0)},B^{(0)}\cap V_\xi,A^{(1)}_\xi,B^{(1)}\cap V_\xi,\hdots,A^{(k-1)}_\xi,B^{(k-1)}\cap V_\xi,A^{(k)}_\xi]\tag*{$(3)$},
	}

	where $\mathcal A^{(i)}=[\vec A^{(i)}]_{\mu_i}\in\ult(\M_i,\mu_i)$ as in the proof of Theorem \ref{theo.ind}.

	\qquad Since $\kappa$ in particular is $\Sigma^1_1$-indescribable, Lemma 3.4.3 of \cite{Abramson} implies that we get a club $C\subset\kappa$ of reflection points of $(2)$. Let $\M_{k+2}\supset\M_{k+1}$ be a weak $\kappa$-model with $\mathcal A^{(k)}\in\M_{k+2}$, where the above $(n+1)$-steps ensured that the $B^{(i)}$'s and the remaining $\mathcal A^{(i)}$'s are all elements of $\M_{k+1}$. In particular, as $C$ is a definable subset in the $\mathcal A^{(i)}$'s and $B^{(i)}$'s we also get that $C\in\M_{k+2}$. Letting $\tilde\mu_{k+2}$ be the associated measure on $\kappa$, $\M_{k+2}$-normality of $\tilde\mu_{k+2}$ ensures that $C\in\tilde\mu_{k+2}$. Now define, for every $\alpha<\kappa$,
	\eq{
		S_\alpha:=\{\xi\in Y\mid\forall i\leq k:\mathcal A^{(i)}\cap V_\alpha=A^{(i)}_\xi\cap V_\alpha\}
	}

	and note that $S_\alpha\in\tilde\mu_{k+2}$ for every $\alpha<\kappa$. Write $\vec S:=\bra{S_\alpha\mid\alpha<\kappa}$ and note that since $\vec S$ is definable it is an element of $\M_{k+2}$ as well. Then $\M_{k+2}$-normality of $\tilde\mu_{k+2}$ ensures that $\triangle\vec S\in\tilde\mu_{k+2}$, so that $C\cap\triangle\vec S\in\tilde\mu_{k+2}$ as well. But letting $\zeta\in C\cap\triangle\vec S$ we see, as in the proof of Theorem \ref{theo.ind}, that
	\eq{
		V_{\zeta+1}\models\varphi[\zeta,A^{(0)}_\zeta,B^{(0)}\cap V_\zeta,A^{(1)}_\zeta,B^{(1)}\cap V_\zeta,\hdots,A^{(k)}_\zeta]
	}

	since $\triangle\vec S\subset Y$, contradicting $(3)$. Hence $X\in\mu$, and since $\M\prec H_\theta$ we have that $\M$ is correct about stationary subsets of $\kappa$, meaning that $\kappa$ is a stationary limit of ${<}\omega$-Ramseys.
}

Now, having established the strength of this large cardinal notion, we move towards complete ineffability. We recall the following definitions.

\defi{
	A collection $R\subset\p(\kappa)$ is a \textbf{stationary class} if
	\begin{enumerate}
		\item $R\neq\emptyset$;
		\item every $A\in R$ is stationary in $\kappa$;
		\item if $A\in R$ and $B\supset A$ then $B\in R$.
	\end{enumerate}
}

\defi{
	A cardinal $\kappa$ is \textbf{completely ineffable} if there is a stationary class $R$ such that for every $A\in R$ and $f:[A]^2\to 2$ there is an $H\in R$ homogeneous for $f$.
}

We then arrive at the following characterisation, influenced by the proof of Theorem 1.3.4 in \cite{Abramson}.

\theo[N.]{
	\label{theo.ineff}
	A cardinal $\kappa$ is completely ineffable if and only if it is coherent ${<}\omega$-Ramsey.
}
\proof{
	$(\Leftarrow)$: Assume $\kappa$ is coherent ${<}\omega$-Ramsey, witnessed by strategies $\bra{\tau_n\mid n<\omega}$. Let $f:[\kappa]^2\to 2$ be arbitrary and form the sequence $\bra{A_\alpha^f\mid\alpha<\kappa}$ as
	\eq{
		A_\alpha^f:=\{\beta>\alpha\mid f(\{\alpha,\beta\})=0\}.
	}

	Let $\M_f$ be a transitive weak $\kappa$-model with $\vec A^f\in\M_f$, and let $\mu_f$ be the associated $\M_f$-measure on $\kappa$ given by $\tau_0$.\footnote{Technically we would have to require that $\M_f\prec H_\theta$ for some regular $\theta>\kappa$ to be able to use $\tau_0$, but note that we could simply get a measure on $\hull^{H_\theta}(\M_f)$ and restrict it to $\M_f$. We will use this throughout the proof.} $1$-Ramseyness of $\kappa$ ensures that $\mu_f$ is normal, meaning $\triangle\mu_f$ is stationary in $\kappa$. Define a new sequence $\vec B^f$ as the $\mu_f$-positive part of $\vec A^f$.\footnote{The \textit{$\mu$-positive part} was defined in Claim \ref{clai.pospart}.} Then $B_\alpha^f\in\mu_f$ for all $\alpha<\kappa$, so that normality of $\mu_f$ implies that $\triangle\vec B^f$ is stationary.
	
	\qquad Let now $\M_f'$ be a new transitive weak $\kappa$-model with $\M_f\subset\M_f'$ and $\mu_f\in\M_f'$, and use $\tau_1$ to get an $\M_f'$-measure $\mu_f'\supset\mu_f$ on $\kappa$. Then $\triangle\vec B^f\cap\{\xi<\kappa\mid A_\xi^f\in\mu_f\}$ and $\triangle\vec B^f\cap\{\xi<\kappa\mid A_\xi^f\notin\mu_f\}$ are both elements of $\M_f'$, so one of them is in $\mu_f'$; set $H_f$ to be that one. Note that $H_f$ is now both stationary in $\kappa$ and homogeneous for $f$.

	\qquad Now let $g:[H_f]^2\to 2$ be arbitrary and again form
	\eq{
		A_\alpha^g:=\{\beta\in H_f\mid\beta>\alpha\land g(\{\alpha,\beta\})=0\}
	}
	
	for $\alpha\in H_f$. Let $\M_{f,g}\supset\M_f'$ be a transitive weak $\kappa$-model with $\vec A^g\in\M_{f,g}$ and use $\tau_2$ to get an $\M_{f,g}$-measure $\mu_{f,g}\supset\mu_f'$ on $\kappa$. As before we then get a stationary $H_{f,g}\in\mu_{f,g}'$ which is homogeneous for $g$. We can continue in this fashion since $\tau_n\subset\tau_{n+1}$ for all $n<\omega$. Define then
	\eq{
		R:=\{A\subset\kappa\mid\exists\vec f:H_{\vec f}\subset A\},
	}

	where the $\vec f$'s range over finite sequences of functions as above; i.e. $f_0:[\kappa]^2\to 2$ and $f_{k+1}:[H_{f_k}]\to 2$ for $k<\omega$.	This is clearly a stationary class which satisfies that whenever $A\in R$ and $g:[A]^2\to 2$, we can find $H\in R$ which is homogeneous for $f$. Indeed, if we let $\vec f$ be such that $H_{\vec f}\subset A$, which exists as $A\in R$, then we can simply let $H:=H_{\vec f,g}$. This shows that $\kappa$ is completely ineffable.

	\qquad $(\Rightarrow)$: Now assume that $\kappa$ is completely ineffable and let $R$ be the corresponding stationary class. We show that $\kappa$ is $n$-Ramsey for all $n<\omega$ by induction, where we inductively make sure that the resulting strategies are coherent as well. Let player I in $\G_0(\kappa)$ play $\M_0$ and enumerate $\p(\kappa)\cap\M_0$ as $\vec A^0\bra{A^0_\alpha\mid\alpha<\kappa}$ such that $A^0_\xi\subset A^0_\zeta$ implies $\xi\leq\zeta$. For $\alpha<\kappa$ define sequences $r_\alpha:\alpha\to 2$ as $r_\alpha(\xi)=1$ iff $\alpha\in A^0_\xi$. Let $<_{\text{lex}}^\alpha$ be the lexicographical ordering on $^\alpha 2$. Define now a colouring $f:[\kappa]^2\to 2$ as
	\eq{
		f(\{\alpha,\beta\}):=\left\{\begin{array}{ll}0 & \text{if }r_{\min(\alpha,\beta)}<_{\text{lex}}^{\min(\alpha,\beta)}r_{\max(\alpha,\beta)}\restr\min(\alpha,\beta)\\ 1 & \text{otherwise}\end{array}\right.
	}

	Let $H_0\in R$ be homogeneous for $f$, using that $\kappa$ is completely ineffable. For $\alpha<\kappa$ consider now the sequence $\bra{r_\xi\restr\alpha\mid\xi\in H_0\land\xi>\alpha}$, which is of length $\kappa$ so there is an $\eta\in[\alpha,\kappa)$ satisfying that $r_\beta\restr\alpha=r_\gamma\restr\alpha$ for every $\beta,\gamma\in H_0$ with $\eta\leq\beta<\gamma$. Define $g:\kappa\to\kappa$ as $g(\alpha)$ being the least such $\eta$, which is then a continuous non-decreasing cofinal function, making the set of fixed points of $g$ club in $\kappa$ -- call this club $C$.

	\qquad Since $H_0$ is stationary we can pick some $\zeta\in C\cap H_0$. As $\zeta\in C$ we get $g(\zeta)=\zeta$, meaning that $r_\beta\restr\zeta=r_\gamma\restr\zeta$ holds for every $\beta,\gamma\in H_0$ with $\zeta\leq\beta<\gamma$. As $\zeta$ is also a member of $H_0$ we can let $\beta:=\zeta$, so that $r_\zeta=r_\gamma\restr\zeta$ holds for every $\gamma\in H_0$, $\gamma>\zeta$.	Now, by definition of $r_\alpha$ we get that for every $\alpha,\gamma\in H_0\cap C$ with $\alpha\leq\gamma$ and $\xi<\alpha$, $\alpha\in A^0_\xi$ iff $\gamma\in A^0_\xi$. Define thus the $\M_0$-measure $\mu_0$ on $\kappa$ as
	\eq{
		\mu_0(A^0_\xi)=1\quad&\text{iff}\quad(\forall\beta\in H_0\cap C)(\beta>\xi\to\beta\in A^0_\xi)\\
		&\text{iff}\quad(\exists\beta\in H_0\cap C)(\beta>\xi\land\beta\in A^0_\xi),
	}

	where the last equivalence is due to the above-mentioned property of $H_0\cap C$. Note that the choice of enumeration implies that $\mu_0$ is indeed a filter. Letting $\vec B=\bra{B_\alpha\mid\alpha<\kappa}$ be the $\mu_0$-positive part of $\vec A^0$, it is also simple to check that $H_0\cap C\subset\triangle\vec B$, making $\mu_0$ normal and hence also both $\M_0$-normal and good, showing that $\kappa$ is $0$-Ramsey.

	\qquad Assume now that $\kappa$ is $n$-Ramsey and let $\bra{\M_0,\mu_0,\hdots,\M_n,\mu_n,\M_{n+1}}$ be a partial play of $\G_{n+1}(\kappa)$. Again enumerate $\p(\kappa)\cap\M_{n+1}$ as $\vec A^{n+1}=\bra{A^{n+1}_\xi\mid\xi<\kappa}$, again satisfying that $\xi\leq\zeta$ whenever $A^{n+1}_\xi\subset A^{n+1}_\zeta$, but also such that given any $\xi<\kappa$ there are $\zeta,\zeta'\in(\xi,\kappa)$ satisfying that $A^{n+1}_\zeta\in\p(\kappa)\cap\M_n$ and $A^{n+1}_{\zeta'}\in(\p(\kappa)\cap\M_{n+1})-\M_n$. The plan now is to do the same thing as before, but we also have to check that the resulting measure extends the previous ones.

	\qquad Let $H_n\in R$ and $C$ be club in $\kappa$ such that $H_n\cap C\subset\triangle\mu_n$, which exist by our inductive assumption. For $\alpha<\kappa$ define $r_\alpha:\alpha\to 2$ as $r_\alpha(\xi)=1$ iff $\alpha\in A^{n+1}_\xi$, and define a colouring $f:[H_n]^2\to 2$ as
	\eq{
		f(\{\alpha,\beta\}):=\left\{\begin{array}{ll}0 & \text{if }r_{\min(\alpha,\beta)}<_{\text{lex}}^{\min(\alpha,\beta)}r_{\max(\alpha,\beta)}\restr\min(\alpha,\beta)\\ 1 & \text{otherwise}\end{array}\right.
	}

	As $H_n\in R$ there is an $H_{n+1}\in R$ homogeneous for $f$. Just as before, define $g:\kappa\to\kappa$ as $g(\alpha)$ being the least $\eta\in[\alpha,\kappa)$ such that $r_\beta\restr\alpha=r_\gamma\restr\alpha$ for every $\beta,\gamma\in H_{n+1}$ with $\eta\leq\beta<\gamma$, and let $D$ be the club of fixed points of $g$. As above we get that given any $\alpha,\gamma\in H_{n+1}\cap D$ with $\alpha\leq\gamma$ and $\xi<\alpha$, $\alpha\in A^{n+1}_\xi$ iff $\gamma\in A^{n+1}_\xi$. Define then the $\M_{n+1}$-measure $\mu_{n+1}$ on $\kappa$ as
	\eq{
		\mu_{n+1}(A^{n+1}_\xi)=1\quad&\text{iff}\quad(\forall\beta\in H_{n+1}\cap D\cap C)(\beta>\xi\to\beta\in A^{n+1}_\xi)\\
		&\text{iff}\quad(\exists\beta\in H_{n+1}\cap D\cap C)(\beta>\xi\land\beta\in A^{n+1}_\xi).
	}

	Then $H_{n+1}\cap D\cap C\subset\triangle\mu_{n+1}$, making $\mu_{n+1}$ normal, $\M_{n+1}$-normal and good, just as before. It remains to show that $\mu_n\subset\mu_{n+1}$. Let thus $A\in\mu_n$ be given, and say $A=A^{n+1}_\xi=A^n_\eta$, where $\vec A^n$ was the enumeration of $\p(\kappa)\cap\M_n$ used at the $n$'th stage. Then by definition of $\mu_n$ we get that for every $\beta\in H_n\cap C$ with $\beta>\eta$, $\beta\in A^n_\eta$. We need to show that
	\eq{
		(\exists\beta\in H_{n+1}\cap D\cap C)(\beta>\xi\land\beta\in A^{n+1}_\xi)
	}

	holds. But here we can simply pick a $\beta>\max(\xi,\eta)$ with $\beta\in H_{n+1}\cap D\cap C\subset H_n\cap C$. This shows that $\mu_n\subset\mu_{n+1}$, making $\kappa$ $(n{+}1)$-Ramsey and thus inductively also coherent ${<}\omega$-Ramsey.
}

\section{The countable case}

This section covers the (strategic) $\gamma$-Ramsey cardinals whenever $\gamma$ has countable cofinality. This case is special because, as mentioned in Section \ref{sect.settingthescene}, we cannot ensure that the final measure is countably complete and so the existence of winning strategies in the $\G_\gamma^\theta(\kappa)$ \textit{might} depend on $\theta$, in contrast with the uncountable cofinality case; see e.g. Question \ref{ques.ctbltheta}.

\subsection*{[Strategic] $\omega$-Ramsey cardinals}

We now move to the strategic $\omega$-Ramsey cardinals and their relationship to the (non-strategic) $\omega$-Ramseys. For this we define a new addition to the family of \textit{virtual cardinals} from \cite{GitmanSchindler}, the \textit{virtually measurable cardinals}.

\defi{
	A cardinal $\kappa$ is \textbf{virtually measurable} if for every regular $\nu>\kappa$ there exists a transitive $M$ and a forcing $\mathbb P$ such that, in $V^{\mathbb P}$, there exists an elementary embedding $j:H_\nu^V\to M$ with $\crit j=\kappa$.
}

We'll need the following well-known lemmata; see Lemma 7.1 in \cite{HolySchlicht} and Lemma 3.1 in \cite{GitmanSchindler} for their proofs.

\qlemm[Ancient Kunen Lemma]{
  \label{lemm.Kunen}
  Let $M\models\zfc^-$ and $j:M\to N$ an elementary embedding with critical point $\kappa$ such that $\kappa+1\subset M\subset N$. Assume that $X\in M$ has $M$-cardinality $\kappa$. Then $j\restr X\in N$.
}

\qlemm[Absoluteness of embeddings on countable structures]{
  \label{lemm.absemb}
  Let $M$ be a countable first-order structure and $j:M\to N$ an elementary embedding. If $W$ is a transitive (set or class) model of (some sufficiently large fragment of) $\zfc$ such that $M$ is countable in $W$ and $N\in W$, then for any finite subset of $M$, $W$ has some elementary embedding $j^*:M\to N$, which agrees with $j$ on that subset. Moreover, if both $M$ and $N$ are transitive $\in$-structures and $j$ has a critical point, we can also assume that $\crit(j^*)=\crit(j)$.
}

\theo[Schindler-N.]{
  \label{theo.virtstrat}
  Every virtually measurable cardinal is strategic $\omega$-Ramsey, and every strategic $\omega$-Ramsey cardinal is virtually measurable in $L$.
}
\proof{
Let $\kappa$ be virtually measurable and fix a regular $\nu>\kappa$, a transitive $M$, a poset $\mathbb P$ and, in $V^{\mathbb P}$, an elementary embedding $\pi:H_\nu^V\to M$ with $\crit\pi=\kappa$. Fix a name $\dot\mu$ and a $\mathbb P$-condition $p$ such that\footnote{Recall that an $M$-measure $\mu$ is \textbf{1-good} if it's weakly amenable and $\ult(M,\mu)$ is well-founded.}
\eq{
  p\forces\godel{\dot\mu\text{ is a 1-good $\check H_\nu$-normal $\check H_\nu$-measure}}
}

We now define a strategy $\sigma$ for player II in $\G_\omega^\nu(\kappa)$ as follows. Whenever player I plays a weak $\kappa$-model $M_n\prec H_\nu^V$, player II fixes $p_n\in\mathbb P$, an $M_n$-measure $\mu_n$ and a function $\pi_n:M_n\to V$ such that $p_0\leq p$, $p_n\leq p_k$ for every $k\leq n$ and that
\eq{
p_n\forces\godel{\dot\mu\cap\check M_n=\check\mu_n\land\check\pi_n=\dot\pi\restr\check M_n}.\tag*{$(1)$}
}

Note that by the Ancient Kunen Lemma \ref{lemm.Kunen} we get that $\pi\restr M_n\in M\subset V$, so such $\pi_n$ always exist in $V$. The $\mu_n$'s also always exist in $V$, by weak amenability of $\mu$. Player II responds to $M_n$ with $\mu_n$. It's clear that the $\mu_n$'s are legal moves for player II, so it remains to show that $\mu_\omega:=\bigcup_{n<\omega}\mu_n$ is good. Assume it's not, so that we have a sequence $\bra{g_n\mid n<\omega}$ of functions $g_n:\kappa\to M_\omega:=\bigcup_{n<\omega}M_n$ such that $g_n\in M_\omega$ and
\eq{
  X_{n+1}:=\{\alpha<\kappa\mid g_{n+1}(\alpha)<g_n(\alpha)\}\in\mu_\omega\tag*{$(2)$}
}

for every $n<\omega$. Without loss of generality we can assume that $g_n,X_n\in M_n$. Then $(2)$ implies that $p_{n+1}\forces\godel{\dot\pi(\check g_{n+1})(\check\kappa)<\dot\pi(\check g_n)(\check\kappa)}$, but by $(1)$ this also means that
\eq{
  p_{n+1}\forces\godel{\check\pi_{n+1}(\check g_{n+1})(\check\kappa)<\check\pi_n(\check g_n)(\check\kappa)},\tag*{$(3)$}
}

so defining, in $V$, the ordinals $\alpha_n:=\pi_n(g_n)(\kappa)$, $(3)$ implies that $\alpha_{n+1}<\alpha_n$ for all $n<\omega$, $\contr$. So $\mu_\omega$ is good, making $\sigma$ a winning strategy and thus also making $\kappa$ strategic $\omega$-Ramsey since $\nu$ was arbitrary.

\qquad Next, let $\kappa$ be strategic $\omega$-Ramsey and fix a winning strategy $\sigma$ for player II in $\G_\omega^\nu(\kappa)$ for a regular $\nu>\kappa$. Let $g\subset\col(\omega,H_\nu^L)$ be $V$-generic and in $V[g]$ fix an elementary chain $\bra{L_{\kappa_n}\mid n<\omega}$ of weak $\kappa$-models $L_{\kappa_n}\prec H_\nu^L$ such that $H_\nu^L\subset\bigcup_{n<\omega}L_{\kappa_n}$, using that $\nu$ is regular and has countable cofinality in $V[g]$. Player II follows $\sigma$, resulting in a $H_\nu^L$-normal $H_\nu^L$-measure $\mu$ on $\kappa$.

\clai{
  $\ult(H_\nu^L,\mu)$ is well-founded.
}

\cproof{
Assume for a contradiction that $\ult(H_\nu^L,\mu)$ is illfounded, witnessed by a sequence $\bra{g_n\mid n<\omega}$ of functions $g_n:\kappa\to\nu$ such that $g_n\in H_\nu^L$ and $\{\alpha<\kappa\mid g_{n+1}(\alpha)<g_n(\alpha)\}\in\mu$. Now, in $V$, define a tree $\T$ of triples $(f,M_f,\mu_f)$ such that $f:\kappa\to\nu$, $M_f$ is a weak $\kappa$-model, $\mu_f$ is an $M_f$-measure on $\kappa$ and letting $f_0<_{\T}\dots <_{\T}f_n=f$ be the $\T$-predecessors of $f$,
\begin{itemize}
  \item $\bra{M_{f_0},\mu_{f_0},\hdots,M_{f_n},\mu_{f_n}}$ is a partial play of $\G_\omega^\nu(\kappa)$ in which player II follows $\sigma$; and
  \item $\{\alpha<\kappa\mid f_{k+1}(\alpha)<f_k(\alpha)\}\in\mu_{k+1}$ for every $k<n$.\\
\end{itemize}

Now, the $g_n$'s induce a cofinal branch through $\T$ in $V[g]$, so by absoluteness of well-foundedness there's a cofinal branch $b$ through $\T$ in $V$ as well. But $b$ now gives us a play of $\G_\omega^\nu(\kappa)$ where player II is following $\sigma$ but player I wins, a contradiction. Thus $\ult(H_\nu^L,\mu)$ is well-founded.
}

Let $j:H_\nu^L\to\ult(H_\nu^L,\mu)\cong M$ be the ultrapower embedding followed by the transitive collapse, so that $M=L_\alpha$ for some $\alpha$ by elementarity. Let now $h\subset\col(\omega,\kappa^{+L})^L$ be $L$-generic, so that $H_\nu^L$ is countable in $L[h]$ and (trivially) $M\in L[h]$. By Lemma \ref{lemm.absemb} we then get that there's an elementary embedding $j^*:H_\nu^L\to M$ in $L[h]$ with critical point $\kappa$. Since we also have that $M\in L$ and as $\nu$ was arbitrary, this makes $\kappa$ virtually measurable in $L$.
}

We get the following immediate corollary.

\qcoro[Schindler-N.]{
	\label{coro.strategic}
  Strategic $\omega$-Ramseys are downwards absolute to $L$, and the existence of a strategic $\omega$-Ramsey cardinal is equiconsistent with the existence of a virtually measurable cardinal. Further, in $L$ the two notions are equivalent.
}

Note also that the proof of Theorem \ref{theo.virtstrat} shows that whenever $\kappa$ is strategic $\omega$-Ramsey then for every regular $\nu>\kappa$ there's a generic extension in which there exists a weakly amenable $H_\nu^V$-normal $H_\nu$-measure on $\kappa$.

\qquad We end this section with a result showing precisely where in the large cardinal hierarchy the strategic $\omega$-Ramsey cardinals and $\omega$-Ramsey cardinals lie, namely that strategic $\omega$-Ramseys are equiconsistent with \textit{remarkables} and $\omega$-Ramseys are strictly below. Theorem 4.8 of \cite{Ramsey2} showed that 2-iterables are limits of remarkables, and our Propositions \ref{prop.tildegame} and \ref{prop.ramseyit} shows that $\omega$-Ramseys are limits of 1-iterables, so that the strategic $\omega$-Ramseys and the $\omega$-Ramseys both lie strictly between the 2-iterables and 1-iterables. It was shown in \cite{HolySchlicht} that $\omega$-Ramseys are consistent with $V=L$. Remarkable cardinals were introduced by \cite{remarkable}, and \cite{GitmanSchindler} showed the following two equivalent formulations.

\defi{
	\label{defi.remarkable}
	A cardinal $\kappa$ is \textbf{remarkable} if one of the two equivalent properties hold:
	\begin{enumerate}
		\item For all $\lambda>\kappa$ there exist $\nu>\lambda$, a transitive set $M$ with $H_\lambda^V\subset M$ and a forcing poset $\mathbb P$, such that in $V^{\mathbb P}$ there's an elementary embedding $\pi:H_\nu^V\to M$ with critical point $\kappa$ and $\pi(\kappa)>\lambda$;
		\item For all $\lambda>\kappa$ there exist $\nu>\lambda$, a transitive set $M$ with ${^\lambda}M\subset M$ and a forcing poset $\mathbb P$, such that in $V^{\mathbb P}$ there's an elementary embedding $\pi:H_\nu^V\to M$ with critical point $\kappa$ and $\pi(\kappa)>\lambda$.
	\end{enumerate}
}

\pagebreak
\theo[N.]{
	\label{theo.remarkable}
	Let $\kappa$ be a virtually measurable cardinal. Then either $\kappa$ is either remarkable in $L$ or $L_\kappa\models\godel{\text{there is a proper class of virtually measurables}}$. In particular, the two notions are equiconsistent.
}
\proof{
	Virtually measurables are downwards absolute to $L$ by Lemma \ref{lemm.absemb}, so we may assume $V=L$. Assume $\kappa$ is not remarkable. This means that there exists some $\lambda>\kappa$ such that for every $\nu>\lambda$, transitive $M$ with $H_\lambda^V\subset M$ and forcing poset $\mathbb P$ it holds that, in $V^{\mathbb P}$, there's no elementary embedding $\pi:H_\nu^V\to M$ with $\crit\pi=\kappa$ and $\pi(\kappa)>\lambda$.
	
	\qquad Fix $\nu:=\lambda^+$ and use that $\kappa$ is virtually $\nu$-measurable to fix a transitive $M$ and a forcing poset $\mathbb P$ such that, in $V^{\mathbb P}$, there's an elementary $\pi:H_\nu^V\to M$. Note that because $M\models V=L$ and $M$ is transitive, $M=L_\alpha$ for some $\alpha\geq\nu$, so that $H_\nu^V=L_\nu\subset M$. This means that $\pi(\kappa)\leq\lambda<\nu$ since we're assuming that $\kappa$ isn't remarkable. Then by restricting the generic embedding to $H_\kappa^V$ we get that $H_\kappa^V\prec H_{\pi(\kappa)}^M=H_{\pi(\kappa)}^V$, using that $\pi(\kappa)<\nu$ and $H_\nu^V=H_\nu^M$ by the above.
	
	\qquad Note that $\pi(\kappa)$ is a cardinal in $H_\nu^V$ since $\pi(\kappa)<\nu$, and as $H_\nu^V\prec_1 V$ we get that $\pi(\kappa)$ \textit{is} a cardinal. But then, again using that $H_{\pi(\kappa)}\prec_1 V$, $\kappa$ is virtually measurable in $H_{\pi(\kappa)}^V$ since being virtually measurable is $\Pi_2$. This means that for every $\xi<\kappa$ it holds that
	\eq{
		H_{\pi(\kappa)}^V\models\exists\alpha>\xi:\godel{\text{$\alpha$ is virtually measurable}},
	}

	implying that $H_\kappa^V\models\godel{\text{There is a proper class of virtually measurables}}$.
}

Now Theorem \ref{theo.remarkable} and Corollary \ref{coro.strategic} yield the following immediate corollary.

\qcoro[Schindler-N.]{
  Let $\kappa$ be strategic $\omega$-Ramsey. Then either $\kappa$ is remarkable in $L$ or otherwise $L_\kappa\models\godel{\text{there is a proper class of strategic $\omega$-Ramseys}}$. In particular, the two notions are equiconsistent.
}

Now, using these results we show that the strategic $\omega$-Ramseys have strictly stronger consistency strength than the $\omega$-Ramseys.

\theo[N.]{
	\label{theo.remlimram}
	Remarkable cardinals are strategic $\omega$-Ramsey limits of $\omega$-Ramsey cardinals.
}
\proof{
	Let $\kappa$ be remarkable. Using property $(ii)$ in the definition of remarkability above we can find a transitive $M$ closed under $2^\kappa$-sequences and a generic elementary embedding $\pi:H_\nu^V\to M$ for some $\nu>2^\kappa$. We will show that $\kappa$ is $\omega$-Ramsey in $M$. Note that remarkables are clearly virtually measurable, and thus by Theorem \ref{theo.virtstrat} also strategic $\omega$-Ramsey; let $\tau_\theta$ be the winning strategy for player II in $\G_\omega^\theta(\kappa)$ for all regular $\theta>\kappa$.

	\qquad In $M$ we fix some regular $\theta>\kappa$ and let $\sigma$ be some strategy for player I in $\G_\omega^\theta(\kappa)^M$. Since $M$ is closed under $2^\kappa$-sequences it means that $\p(\p(\kappa))\subset M$ and thus that $M$ contains all possible filters on $\kappa$. We let player II follow $\tau$, which produces a play $\sigma*\tau$ in which player II wins. But all player II's moves are in $\p(\p(\kappa))$ and hence in $M$, and as $M$ is furthermore closed under $\omega$-sequences, $\sigma*\tau\in M$. This means that $M$ sees that $\sigma$ is not winning, so $\kappa$ is $\omega$-Ramsey in $M$.
	
	\qquad This also implies that $\kappa$ is a limit of $\omega$-Ramseys in $H_\nu$. But as $\kappa$ is remarkable it holds that $H_\kappa\prec_2 V$, in analogy with the same property for strongs and supercompacts, and as being $\omega$-Ramsey is a $\Pi_2$-notion this means that $\kappa$ \textit{is} a limit of $\omega$-Ramseys.
}

This immediately yields the following corollary.

\coro[Schindler-N.]{
	If $\kappa$ is a strategic $\omega$-Ramsey cardinal then
	\eq{
		L_\kappa\models\godel{\text{there is a proper class of $\omega$-Ramseys}}.\tag*{$\dashv$}
	}
}

\subsection*{$(\omega,\alpha)$-Ramsey cardinals}

A natural generalisation of the $\gamma$-Ramsey definition is to require more iterability of the last measure. Of course, by Proposition \ref{prop.tildegame} we have that $\G_\gamma(\kappa,\zeta)$ is equivalent to $\G_\gamma(\kappa)$ when $\cof\gamma>\omega$ so the next definition is only interesting whenever $\cof\gamma=\omega$.

\defi[N.]{
	Let $\alpha,\beta$ be ordinals. Then a cardinal $\kappa$ is \textbf{$(\alpha,\beta)$-Ramsey} if player I does not have a winning strategy in $\G_\alpha^\theta(\kappa,\beta)$ for all regular $\theta>\kappa$.\footnote{Note that an $\alpha$-Ramsey cardinal is the same as an $(\alpha,0)$-Ramsey cardinal.}
}

\defi[Gitman]{
	A cardinal $\kappa$ is \textbf{$\alpha$-iterable} if for every $A\subset\kappa$ there exists a \textit{transitive} weak $\kappa$-model $\M$ with $A\in\M$ and an $\alpha$-good $\M$-measure $\mu$ on $\M$.
}

\prop{
	\label{prop.ramseyit}
	If $\beta>0$ then every $(\alpha,\beta)$-Ramsey is a $\beta$-iterable stationary limit of $\beta$-iterables.
}
\proof{
	Let $(\M,\in,\mu)$ be a result of a play of $\G_\alpha^{\kappa^+}(\kappa,\beta)$ in which player II won. Then the transitive collapse of $(\M,\in,\mu)$ witnesses that $\kappa$ is $\beta$-iterable, since $\mu$ is $\beta$-good by definition of $\G_\alpha^{\kappa^+}(\kappa,\beta)$.
	
	\qquad That $\kappa$ is $\beta$-iterable is reflected to some $H_\theta$, so let now $(\N,\in,\nu)$ be a result of a play of $\G_\alpha^\theta(\kappa,\beta)$ in which player II won. Then $\N\prec H_\theta$, so that $\kappa$ is also $\beta$-iterable in $\N$. Since being $\beta$-iterable is witnessed by a subset of $\kappa$ and $\beta>0$ implies\footnote{Recall that $\beta$-good for $\beta>0$ in particular implies weak amenability.} that we get a $\kappa$-powerset preserving $j:\N\to\P$, $\P$ also thinks that $\kappa$ is $\beta$-iterable, making $\kappa$ a stationary limit of $\beta$-iterables by elementarity.
}

We now move towards Theorem \ref{theo.upperlimit} which gives an upper consistency bound for the $(\omega,\alpha)$-Ramseys. We first recall a few definitions and a folklore lemma.

\defi{
	For an infinite ordinal $\alpha$, a cardinal $\kappa$ is \textbf{$\alpha$-Erd\H os} for $\alpha\leq\kappa$ if given any club $C\subset\kappa$ and regressive $c:[C]^{<\omega}\to\kappa$ there is a set $H\in[C]^\alpha$ homogeneous for $c$; i.e. that $\abs{c``[H]^n}\leq 1$ holds for every $n<\omega$.
}

\defi{
	\label{defi.remind}
	A set of indiscernibles $I$ for a structure $\M=(M,\in,A)$ is \textbf{remarkable} if $I-\iota$ is a set of indiscernibles for $(M,\in,A,\bra{\xi\mid\xi<\iota})$ for every $\iota\in I$.
}

\lemm[Folklore]{
	\label{lemma}
	Let $\kappa$ be $\alpha$-Erd\H os where $\alpha\in[\omega,\kappa]$ and let $C\subset\kappa$ be club. Then any structure $\M$ in a countable language $\mathcal L$ with $\kappa+1\subset\M$ has a remarkable set of indiscernibles $I\in[C]^\alpha$.
}
\proof{
	Let $\bra{\varphi_n\mid n<\omega}$ enumerate all $\mathcal L$-formulas and define $c:[C]^{<\omega}\to\kappa$ as follows. For an increasing sequence $\alpha_1<\cdots<\alpha_{2n}\in C$ let
	\eq{
		c(\{\alpha_1,\hdots,\alpha_{2n}\}):=&\text{ the least }\lambda<\alpha_1\text{ such that }\exists\delta_1<\cdots\delta_k\exists m<\omega:\lambda=\bra{m,\delta_1,\hdots,\delta_k}\land\\
		&\M\not\models\varphi_m[\vec\delta,\alpha_1,\hdots,\alpha_n]\leftrightarrow\varphi_m[\vec\delta,\alpha_{n+1},\hdots,\alpha_{2n}]
	}

	if such a $\lambda$ exists, and $c(s)=0$ otherwise. Clearly $c$ is regressive, so since $\kappa$ is $\alpha$-Erd\H os we get a homogeneous $I\in[C]^\alpha$ for $c$; i.e. that $\abs{c``[I]^n}\leq 1$ for every $n<\omega$. Then $c(\{\alpha_1,\hdots,\alpha_{2n}\})=0$ for every $\alpha_1,\hdots,\alpha_{2n}\in I$, as otherwise there exists an $m<\omega$ and $\delta_1<\cdots\delta_k$ such that for any $\alpha_1<\hdots<\alpha_{2n}\in I$,
	\eq{
		\M\not\models\varphi_m[\vec\delta,\alpha_1,\hdots,\alpha_n]\leftrightarrow\varphi_m[\vec\delta,\alpha_{n+1},\hdots,\alpha_{2n}].\tag*{$(\dagger)$}
	}
	
	But then simply pick $\alpha_1<\hdots\alpha_{2n}<\alpha_1'<\cdots<\alpha_{2n}'$ so that both $\{\alpha_1,\hdots,\alpha_{2n}\}$ and $\{\alpha_1',\hdots,\alpha_{2n}'\}$ witnesses $(\dagger)$; then either $\{\alpha_1,\hdots,\alpha_n,\alpha_1',\alpha_n'\}$ or $\{\alpha_1,\hdots,\alpha_n,\alpha_{n+1}',\hdots,\alpha_{2n}'\}$ also witnesses that $(\dagger)$ fails, $\contr$.
}

\theo[N.]{
	\label{theo.upperlimit}
	Let $\alpha\in[\omega,\omega_1]$ be additively closed. Then any $\alpha$-Erd\H os cardinal is a limit of $(\omega,\alpha)$-Ramsey cardinals.
}
\proof{
	Let $\kappa$ be $\alpha$-Erd\H os, $\theta>\kappa$ a regular cardinal and $\beta<\kappa$ any ordinal. Use the above Lemma \ref{lemma} to get a set of remarkable indiscernibles $I\in[\kappa]^\alpha$ for the structure $(H_\theta,\in,\bra{\xi\mid\xi<\beta})$, and let $\iota\in I$ be the least indiscernible in $I$. We will show that player I has no winning strategy in $\G_\omega^\theta(\iota,\alpha)$, so by the proof of Theorem 5.5(d) in \cite{HolySchlicht} it suffices to find a weak $\iota$-model $\M\prec H_\theta$ and an $\alpha$-good $\M$-measure on $\iota$. Define
	\eq{
		\M:=\hull^{H_\theta}(\iota\cup I)\prec H_\theta
	}
	
	and let $\pi:I\to I$ be the right-shift map. Since $I$ is remarkable, $I$ ($=I-\iota$) is a set of indiscernibles for the structure $(H_\theta,\in,\bra{\xi\mid\xi<\iota})$, so that $\pi$ induces an elementary embedding $j:\M\to\M$ with $\crit j=\iota$, given as
	\eq{
		j(\tau^{\M}[\vec\xi,\iota_{i_0},\hdots,\iota_{i_k}]):=\tau^{\M}[\vec\xi,\iota_{i_0+1},\hdots,\iota_{i_k+1}],
	}

	with $\vec\xi\subset\iota$. Since $j$ is trivially $\iota$-powerset preserving we get that $\M\prec H_\theta$ is a weak $\iota$-model satisfying $\zfc^-$ with a 1-good $\M$-measure $\mu_j$ on $\iota$. Furthermore, as we can linearly iterate $\M$ simply by applying $j$ we get an $\alpha$-iteration of $\M$ since there are $\alpha$-many indiscernibles. Note that at limit stages $\gamma<\alpha$ our iteration sends $\tau^{\M}[\vec\xi,\iota_{i_0},\hdots,\iota_{i_k}]$ to $\tau^{\M}[\vec\xi,\iota_{i_0+\gamma},\hdots,\iota_{i_k+\gamma}]$ so here we are using that $\alpha$ is additively closed.
	
	\qquad This shows that player I has no winning strategy in $\G_\omega^\theta(\iota,\alpha)$. Since $\iota>\beta$ and $\beta<\kappa$ was arbitrary, $\kappa$ is a limit of $\eta$ such that player I has no winning strategy in $\G_\omega^\theta(\eta,\alpha)$. If we repeat this procedure for all regular $\theta>\kappa$ we get by the pidgeon hole principle that $\kappa$ is a limit of $(\omega,\alpha)$-Ramsey cardinals.
}

As Theorem 4.5 in \cite{GitmanSchindler} shows that $(\alpha{+}1)$-iterable cardinals have $\alpha$-Erd\H os cardinals below them for $\alpha\geq\omega$ additively closed, this shows that the $(\omega,\alpha)$-Ramseys form a strict hierarchy. Further, as $\alpha$-Erd\H os cardinals are consistent with $V=L$ when $\alpha<\omega_1^L$ and $\omega_1$-iterable cardinals aren't consistent with $V=L$, we also get that $(\omega,\alpha)$-Ramsey cardinals are consistent with $V=L$ if $\alpha<\omega_1^L$ and that they aren't if $\alpha=\omega_1$.

\subsection*{[Strategic] $(\omega{+}1)$-Ramsey cardinals}

The next step is then to consider $(\omega{+}1)$-Ramseys, which turn out to cause a considerable jump in consistency strength. We first need the following result which is implicit in \cite{Mitchell} and in the proof of Lemma 1.3 in \cite{Kapplications} --- see also \cite{Dodd} and \cite{Ramsey1}.

\qtheo[Dodd, Mitchell]{
	\label{theo.ramseycondition}
	A cardinal $\kappa$ is Ramsey if and only if every $A\subset\kappa$ is an element of a weak $\kappa$-model $\M$ such that there exists a weakly amenable countably complete $\M$-measure on $\kappa$.
}

The following theorem then supplies us with a lower bound for the strength of the $(\omega{+}1)$-Ramsey cardinals. It should be noted that a better lower bound will be shown in Theorem \ref{theo.plusone}, but we include this Ramsey lower bound as well for completeness.

\theo[N.]{
	\label{theo.ramlimram}
	Every $(\omega{+}1)$-Ramsey cardinal is a Ramsey limit of Ramseys.
}
\proof{
	Let $\kappa$ be $(\omega{+}1)$-Ramsey and $A\subset\kappa$. Let $\sigma$ be a strategy for player I in $\G_{\omega+1}^{\kappa^+}(\kappa)$ satisfying that whenever $\vec\M_\alpha*\vec\mu_\alpha$ is consistent with $\sigma$ it holds that $A\in\M_0$ and $\mu_\alpha\in\M_{\alpha+1}$ for all $\alpha\leq\omega$. Then $\sigma$ isn't winning as $\kappa$ is $(\omega{+}1)$-Ramsey, so we may fix a play $\sigma*\vec\mu_\alpha$ of $\G_{\omega+1}^{\kappa^+}(\kappa)$ in which player II wins. Then by the choice of $\sigma$ we get that $\mu_\omega$ is a weakly amenable $\M_\omega$-measure on $\kappa$, and by the rules of $\G_{\omega+1}^{\kappa^+}(\kappa)$ it's also countably complete (it's even normal), which makes $\kappa$ Ramsey by the above Theorem \ref{theo.ramseycondition}.

	\qquad Since $\kappa$ is Ramsey, $\M_\omega\models\godel{\kappa\text{ is Ramsey}}$ as well. Letting $j:\M_\omega\to\N$ be the $\kappa$-powerset preservering embedding induced by $\mu_\omega$, we also get that $\N\models\godel{\kappa\text{ is Ramsey}}$ by $\kappa$-powerset preservation. This then implies that $\kappa$ \textit{is} a stationary limit of Ramsey cardinals inside $\M_\omega$, and thus also in $V$ by elementarity.
}

As for the \textit{consistency} strength of the strategic $(\omega{+}1)$-Ramsey cardinals, we get the following result that they reach a measurable cardinal. The proof of the following is closely related to the proof due to Silver and Solovay that player II having a winning strategy in the \textit{cut and choose game} is equiconsistent with a measurable cardinal --- see e.g. p. 249 in \cite{c&c}.

\theo[N.]{
	\label{theo.omegaplusone}
	If $\kappa$ is a strategic $(\omega{+}1)$-Ramsey cardinal then, in $V^{\col(\omega,2^\kappa)}$, there's a transitive class $N$ and an elementary embedding $j:V\to N$ with $\crit j=\kappa$. In particular, the existence of a strategic $(\omega{+}1)$-Ramsey cardinal is equiconsistent with the existence of a measurable cardinal.
}
\proof{
	Set $\mathbb P:=\col(\omega,2^\kappa)$ and let $\sigma$ be player II's winning strategy in $\G_{\omega+1}^{\kappa^+}(\kappa)$. Let $\dot\M$ be a $\mathbb P$-name of an $\omega$-sequence $\bra{\M_n\mid n<\omega}$ of weak $\kappa$-models $\M_n\in V$ such that $\M_n\prec H_{\kappa^+}^V$ and $\p(\kappa)^V\subset\bigcup_{n<\omega}\M_n$, and let $\dot\mu$ be a $\mathbb P$-name for the $\omega$-sequence of $\sigma$-responses to the $\M_n$'s in $\G_{\omega+1}^{\kappa^+}(\kappa)^V$.

	\qquad Assume that there's a $\mathbb P$-condition $p$ which forces the generic ultrapower $\ult(V,\bigcup_n\dot\mu_n)$ to be illfounded, meaning that we can fix a $\mathbb P$-name $\dot f$ for an $\omega$-sequence $\bra{f_n\mid n<\omega}$ such that
	\eq{
		p\forces\dot X_n:=\{\alpha<\kappa\mid\dot f_{n+1}(\alpha)<\dot f_n(\alpha)\}\in\bigcup_{n<\omega}\dot\mu_n.
	}

	Now, in $V$, we fix some large regular $\theta\gg\kappa$ and a countable $\N\prec H_\theta$ such that $\dot\M,\dot\mu,\dot f,H_{\kappa^+}^V,\sigma,p\in\N$. We can find an $\N$-generic $g\subset\mathbb P^{\N}$ in $V$ with $p\in g$ since $\N$ is countable, so that $\N[g]\in V$. But the play $\dot\M^g_n*\dot\mu^g_n$ is a play of $\G_\omega^{\kappa^+}(\kappa)^V$ which is according to $\sigma$, meaning that $\bigcup_{n<\omega}\dot\mu^g_n$ is normal and in particular countably complete (in $V$). Then $\bigcap_{n<\omega}\dot X_n^g\neq\emptyset$, but if $\alpha\in\bigcap_{n<\omega}\dot X_n^g$ then $\bra{\dot f^g_n(\alpha)\mid n<\omega}$ is a strictly decreasing $\omega$-sequence of ordinals, $\contr$. This means that $\ult(V,\bigcup_n\mu_n)$ is indeed wellfounded.

	\qquad This conclusion is well-known to imply that $\kappa$ is a measurable in an inner model; see e.g. Lemma 4.2 in \cite{KellnerShelah}.
}

The above Theorem \ref{theo.omegaplusone} then answers Question 9.2 in \cite{HolySchlicht} in the negative, asking if $\lambda$-Ramseys are strategic $\lambda$-Ramseys for uncountable cardinals $\lambda$, as well as answering Question 9.7 from the same paper in the positive, asking whether strategic fully Ramseys are equiconsistent with a measurable.

\section{The general case}

\subsection*{Gitman's cardinals}

In this subsection we define the strongly- and super Ramsey cardinals from \cite{Ramsey1} and investigate further connections between these and the $\alpha$-Ramsey cardinals. First, a definition.

\defi[Gitman]{
	A cardinal $\kappa$ is \textbf{strongly Ramsey} if every $A\subset\kappa$ is an element of a transitive $\kappa$-model $\M$ with a weakly amenable $\M$-normal $\M$-measure $\mu$ on $\kappa$. If furthermore $\M\prec H_{\kappa^+}$ then we say that $\kappa$ is \textbf{super Ramsey}.
}

Note that since the model $\M$ in question is a \textit{$\kappa$-model} it is closed under countable sequences, so that the measure $\mu$ is automatically countably complete. The definition of the strongly Ramseys is thus exactly the same as the characterisation of Ramsey cardinals, with the added condition that the model is closed under ${<}\kappa$-sequences. \cite{Ramsey1} shows that every super Ramsey cardinal is a strongly Ramsey limit of strongly Ramsey cardinals, and that $\kappa$ is strongly Ramsey iff every $A\subset\kappa$ is an element of a transitive $\kappa$-model $\M\models\zfc$ with a weakly amenable $\M$-normal $\M$-measure $\mu$ on $\kappa$.

\qquad Now, a first connection between the $\alpha$-Ramseys and the strongly- and super Ramseys is the result in \cite{HolySchlicht} that fully Ramsey cardinals are super Ramsey limits of super Ramseys. The following result then shows that the strongly- and super Ramseys are sandwiched between the almost fully Ramseys and the fully Ramseys.

\theo[N.-W.]{
  Every strongly Ramsey cardinal is a stationary limit of almost fully Ramseys.
}
\proof{
  Let $\kappa$ be strongly Ramsey and let $\M\models\zfc$ be a transitive $\kappa$-model with $V_\kappa\in\M$ and $\mu$ a weakly amenable $\M$-normal $\M$-measure. Let $\gamma<\kappa$ have uncountable cofinality and $\sigma\in\M$ a strategy for player I in $\G_\gamma(\kappa)^{\M}$. Now, whenever player I plays $\M_\alpha\in\M$ let player II play $\mu\cap\M_\alpha$, which is an element of $\M$ by weak amenability of $\mu$. As $\M^{<\kappa}\subset\M$ the resulting play is inside $\M$, so $\M$ sees that $\sigma$ is not winning.

	\qquad Now, letting $j_\mu:\M\to\N$ be the induced embedding, $\kappa$-powerset preservation of $j_\mu$ implies that $\mu$ is also a weakly amenable $\N$-normal $\N$-measure on $\kappa$. This means that we can copy the above argument to ensure that $\kappa$ is also almost fully Ramsey in $\N$, entailing that it is a stationary limit of almost fully Ramseys in $\M$. But note now that $\lambda$ is almost fully Ramsey iff it is almost fully Ramsey in a transitive $\zfc$-model containing $H_{(2^\lambda)^+}$ as an element by Theorem 5.5(e) in \cite{HolySchlicht}, so that $\kappa$ being inaccessible, $V_\kappa\in\M$ and $\M$ being transitive implies that $\kappa$ really \textit{is} a stationary limit of almost fully Ramseys.
}

\subsection*{Downwards absoluteness to $K$}

Lastly, we consider the question of whether the $\alpha$-Ramseys are downwards absolute to $K$, which turns out to at least be true in many cases. The below Theorem \ref{theo.downK} then also answers Question 9.4 from \cite{HolySchlicht} in the positive, asking whether $\alpha$-Ramseys are downwards absolute to the Dodd-Jensen core model for $\alpha\in[\omega,\kappa]$ a cardinal. We first recall the definition of $0^\pistol$.

\defi{
	$0^\pistol$ is ``the sharp for a strong cardinal'', meaning the minimal sound active mouse $\M$ with $\M\l\crit(\dot F^{\M})\models\godel{\text{There exists a strong cardinal}}$, with $\dot F^{\M}$ being the top extender of $\M$.
}

\theo[N.-W.]{
	\label{theo.downK}
	Assume $0^\pistol$ does not exist. Let $\lambda$ be a limit ordinal with uncountable cofinality and let $\kappa$ be $\lambda$-Ramsey. Then $K\models\godel{\kappa\text{ is a $\lambda$-Ramsey cardinal}}$.
}
\proof{
	Note first that $\kappa^{+K}=\kappa^+$ by \cite{Schindler}, since $\kappa$ in particular is weakly compact.	Let $\sigma\in K$ be a strategy for player I in $\G_\lambda^{\kappa^+}(\kappa)^K$, so that a play following $\sigma$ will produce weak $\kappa$-models $\M\prec K\l\kappa^+$. We can then define a strategy $\tilde\sigma$ for player I in $\G_\lambda^{\kappa^+}(\kappa)$ as follows. Firstly let $\tilde\sigma(\emptyset):=\hull^{H_{\kappa^+}}(K\l\kappa\cup\sigma(\emptyset))$. Assuming now that $\bra{\tilde\M_\alpha,\tilde\mu_\alpha\mid\alpha<\gamma}$ is a partial play of $\G_\lambda^{\kappa^+}(\kappa)$ which is consistent with $\tilde\sigma$, we have two cases. If $\tilde\mu_\alpha\in K$ for every $\alpha<\gamma$ then let $\bra{\M_\alpha\mid\alpha<\gamma}$ be the corresponding models played in $\G_\lambda^{\kappa^+}(\kappa)^K$ from which the $\tilde\M_\alpha$'s are derived and let
	\eq{
		\tilde\sigma(\bra{\tilde\M_\alpha,\tilde\mu_\alpha\mid\alpha<\gamma}):=\hull^{H_{\kappa^+}}(K\l\kappa\cup\sigma(\bra{\M_\alpha,\tilde\mu_\alpha\mid\alpha<\gamma})),
	}

	and otherwise let $\tilde\sigma$ play arbitrarily. As $\kappa$ is $\lambda$-Ramsey (in $V$) there exists a play $\bra{\tilde\M_\alpha,\tilde\mu_\alpha\mid\alpha\leq\lambda}$ of $\G_\lambda^{\kappa^+}(\kappa)$ which is consistent with $\tilde\sigma$ in which player II won. Note that $\tilde\M_\lambda\cap K\l\kappa^+\prec K\l\kappa^+$ so let $\N$ be the transitive collapse of $\tilde\M_\lambda\cap K\l\kappa^+$. But if $j:\N\to K\l\kappa^+$ is the uncollapse then $\crit j$ is both an $\N$-cardinal and also $>\kappa$ because we ensured that $K\l\kappa\subset\N$. This means that $j=\id$ because $\kappa$ is the largest $\N$-cardinal by elementarity in $K\l\kappa^+$, so that $\tilde\M_\lambda\cap K\l\kappa^+=\N$ is a transitive elementary substructure of $K\l\kappa^+$, making it an initial segment of $K$.
	
	\qquad Now, since $\mu:=\tilde\mu_\lambda$ is a countably complete weakly amenable $K\l o(\N)$-measure\footnote{Here we use that $\N\pinit K$.}, the ``beaver argument''\footnote{See Lemmata 7.3.7--7.3.9 and 8.3.4 in \cite{Zeman} for this argument.} shows that $\mu\in K$, so that we can then define a strategy $\tau$ for player II in $\G_\lambda^{\kappa^+}(\kappa)^K$ as simply playing $\mu\cap\N\in K$ whenever player I plays $\N$. Since $\mu=\tilde\mu_\lambda$ we also have that $\mu\cap\M_\alpha=\tilde\mu_\alpha\cap\M_\alpha$, so that $\sigma$ will eventually play $\N$, making $\tau$ win against $\sigma$.\footnote{Note that $\tau$ is not necessarily a winning strategy --- all we know is that it is winning against this particular strategy $\sigma$.}
}

Note that the only thing we used $\cof\lambda>\omega$ for in the above proof was to ensure that $\mu$ was countably complete. If now $\kappa$ instead was either genuine- or normal $\alpha$-Ramsey for any limit ordinal $\alpha$ then $\mu_\alpha$ would also be countably complete and weakly amenable, so the same proof shows the following.

\qcoro[N.-W.]{
	\label{coro.downK}
	Assume $0^\pistol$ does not exist and let $\alpha$ be any limit ordinal. Then every genuine- and every normal $\alpha$-Ramsey cardinal is downwards absolute to $K$. In particular, if $\alpha$ is a limit of limit ordinals then every ${<}\alpha$-Ramsey cardinal is downwards absolute to $K$ as well.
}

\subsection*{Indiscernible games}

We now move to the strategic versions of the $\alpha$-Ramsey hierarchy. The first thing we want to do is define \textit{$\alpha$-very Ramsey cardinals}, introduced in \cite{SharpeWelch}, and show the tight connection between these and the strategic $\alpha$-Ramseys. We need a few more definitions. Recall the definition of a remarkable set of indiscernibles from Definition \ref{defi.remind}.

\defi{
	A \textbf{good set of indiscernibles} for a structure $\M$ is a set $I\subset\M$ of remarkable indiscernibles for $\M$ such that $\M\l\iota\prec\M$ for any $\iota\in I$.
}

\defi[Sharpe-W.]{
	Define the \textbf{indiscernible game} $G^I_\gamma(\kappa)$ in $\gamma$ many rounds as follows
	\game{\M_0}{I_0}{\M_1}{I_1}{\M_2}{I_2}{\cdots}{\cdots}

	Here $\M_\alpha$ is an amenable structure of the form $(J_\kappa[A],\in,A)$ for some $A\subset\kappa$, $I_\alpha\in[\kappa]^\kappa$ is a good set of indiscernibles for $\M_\alpha$ and the $I_\alpha$'s are $\subset$-decreasing. Player II wins iff they can continue playing through all the rounds.
}

\defi[Sharpe-W.]{
	A cardinal $\kappa$ is \textbf{$\gamma$-very Ramsey} if player II has a winning strategy in the game $G^I_\gamma(\kappa)$.
}

The next couple of results concerns the connection between the strategic $\alpha$-Ramseys and the $\alpha$-very Ramseys. We start with the following.

\theo[N.]{
	\label{theo.plusone}
	Every $(\omega{+}1)$-Ramsey is an $\omega$-very Ramsey stationary limit of $\omega$-very Ramseys.
}
\proof{
	Let $\kappa$ be $(\omega{+}1)$-Ramsey. We will describe a winning strategy for player II in the indiscernible game $G_\omega^I(\kappa)$. If player I plays $\M_0=(J_\kappa[A_0],\in,A_0)$ in $G_\omega^I(\kappa)$ then let player I in $\G_{\omega+1}^{\kappa^+}(\kappa)$ play
	\eq{
		\h_0:=\hull^{H_{\kappa^+}}(J_\kappa[A_0]\cup\{\M_0,\kappa,A_0\})\prec H_{\kappa^+}.
	}

	Let player I now follow a strategy in $\G_{\omega+1}^{\kappa^+}(\kappa)$ which starts off with $\h_0$ and ensures that, whenever $\vec\M_\alpha*\vec\mu_\alpha$ is consistent with player I's strategy, then $\mu_\alpha\in\M_{\alpha+1}$ for all $\alpha\leq\omega$. Since player II is not losing in $\G_{\omega+1}^{\kappa^+}(\kappa)$ there is a play $\vec\M_\alpha*\vec\mu_\alpha$ in which player I follows this strategy just described and where player II wins -- write $\h_0^{(\alpha)}:=\M_\alpha$ and $\mu_0^{(\alpha)}:=\mu_\alpha$ for the models and measures in this play.
	\game{\h_0^{(0)}}{\mu_0^{(0)}}{\cdots}{\cdots}{\h_0^{(\omega)}}{\mu_0^{(\omega)}}{\h_0^{(\omega+1)}}{\mu_0^{(\omega+1)}}

	By the choice of player I's strategy we get that $\mu_0^{(\omega)}$ is both weakly amenable, and it's also countably complete by the rules of $\G_{\omega+1}^{\kappa^+}(\kappa)$ (it's even normal). Now Lemma 2.9 of \cite{SharpeWelch} gives us a set of good indiscernibles $I_0\in\mu_0^{(\omega)}$ for $\M_0$, as $\M_0\in\h_0^{(\omega)}$ and $\mu_0^{(\omega)}$ is a countably complete weakly amenable $\h_0^{(\omega)}$-normal $\h_0^{(\omega)}$-measure on $\kappa$. Let player II play $I_0$ in $G_\omega^I(\kappa)$. Let now $\M_1=(J_\kappa[A_1],\in,A_1)$ be the next play by player I in $G_\omega^I(\kappa)$.
	\game{\M_0}{I_0}{\M_1}{}{}{}{}{}

	Since $\mu_0^{(\omega)}=\bigcup_n\mu_0^{(n)}$ we must have that $I_0\in\mu_0^{(n_0)}$ for some $n_0<\omega$. In the $(n_0{+}1)$'st round of $\G_{\omega+1}^{\kappa^+}(\kappa)$ we change player I's strategy and let player I play
	\eq{
		\h_1:=\hull^{H_{\kappa^+}}(J_\kappa[A_0]\cup\{\M_0,\M_1,\kappa,A_0,A_1,\bra{\h_0^{(k)},\mu_0^{(k)}\mid k\leq n_0}\})\prec H_{\kappa^+}
	}

	and otherwise continues following some strategy, as long as the measures played by player II keep being elements of the following models.	Our play of the game $\G_{\omega+1}^{\kappa^+}(\kappa)$ thus looks like the following so far.

	\game{\h_0^{(0)}}{\mu_0^{(0)}}{\cdots}{\cdots}{\h_0^{(n_0)}}{\mu_0^{(n_0)}}{\h_1}{}

	Now player II in $\G_{\omega+1}^{\kappa^+}(\kappa)$ is not losing at round $n_0$, so there is a play extending the above in which player I follows their revised strategy and in which player II wins. As before we get a set $I_1'\in\mu_1^{(n_1)}$ of good indiscernibles for $\M_1$, where $n_1<\omega$. Since $I_0\in\mu_0^{(n_0)}\subset\mu_1^{(n_1)}$ we can let player II in $G_\omega^I(\kappa)$ play $I_1:=I_0\cap I_1'\in\mu_1^{(n_1)}$. Continuing like this, player II can keep playing throughout all $\omega$ rounds of $G_\omega^I(\kappa)$, making $\kappa$ $\omega$-very Ramsey.

	\qquad As for showing that $\kappa$ is a stationary limit of $\omega$-very Ramseys, let $\M\prec H_{\kappa^+}$ be a weak $\kappa$-model with a weakly amenable countably complete $\M$-normal $\M$-measure $\mu$ on $\kappa$, which exists by Theorem \ref{theo.ramlimram} as $\kappa$ is $(\omega{+}1)$-Ramsey. Then by elementarity $\M\models\godel{\kappa\text{ is $\omega$-very Ramsey}}$ and since $\kappa$ being $\omega$-very Ramsey is absolute between structures having the same subsets of $\kappa$ it also holds in the $\mu$-ultrapower, meaning that $\kappa$ is a stationary limit of $\omega$-very Ramseys by elementarity. 
}

The above proof technique can be generalised to the following.

\theo[N.]{
	\label{theo.stratvery}
	For limit ordinals $\alpha$, every coherent ${<}\omega\alpha$-Ramsey is $\omega\alpha$-very Ramsey.
}
\proof{
	This is basically the same proof as the proof of Theorem \ref{theo.plusone}. We do the ``going-back'' trick in $\omega$-chunks, and at limit stages we continue our non-losing strategy in $\G_{\omega\alpha}^{\kappa^+}(\kappa)$ by using our winning strategy, which we have available as we are assuming coherent ${<}\omega\alpha$-Ramseyness. We need $\alpha$ to be a limit ordinal for this to work, as otherwise we would be in trouble in the last $\omega$-chunk, as we cannot just extend the play to get a countably complete measure, which we need to use the proof of Theorem \ref{theo.plusone}.
}

As for going from the $\alpha$-very Ramseys to the strategic $\alpha$-Ramseys we got the following.

\theo[N.]{
	\label{theo.succverystrat}
	For $\gamma$ any ordinal, every coherent ${<}\gamma$-very Ramsey\footnote{Here the coherency again just means that the winning strategies $\sigma_\alpha$ for player II in $G_\alpha^I(\kappa)$ are $\subset$-increasing.} is coherent ${<}\gamma$-Ramsey.\footnote{Here a ``coherent ${<}\gamma$-very Ramsey cardinal'' is defined from $\gamma$-very Ramseys in the same way as coherent ${<}\gamma$-Ramsey cardinals is defined from $\gamma$-Ramseys. When $\gamma$ is a limit ordinal then coherent ${<}\gamma$-very Ramseys are precisely the same as $\gamma$-very Ramseys, so this is solely to ``subtract one'' when $\gamma$ is a successor ordinal --- i.e. a coherent ${<}(\gamma+1)$-very Ramsey cardinal is the same thing as a $\gamma$-very Ramsey cardinal.}
}
\proof{
	The reason why we work with ${<}\gamma$-Ramseys here is to ensure that player II only has to satisfy a closed game condition (i.e. to continue playing throughout all the rounds). If $\gamma=\beta+1$ then set $\zeta:=\beta$ and otherwise let $\zeta:=\gamma$. Let $\kappa$ be $\zeta$-very Ramsey and let $\tau$ be a winning strategy for player II in $G_\zeta^I(\kappa)$. Let $\M_\alpha\prec H_\theta$ be any move by player I in the $\alpha$'th round of $\G_\zeta(\kappa)$. Let $A_\alpha\subset\kappa$ encode all subsets of $\kappa$ in $\M_\alpha$ and form now
	\eq{
		\N_\alpha:=(J_\kappa[A_\alpha],\in,A_\alpha),
	}

	which is a legal move for player I in $G_\zeta^I(\kappa)$, yielding a good set of indiscernibles $I_\alpha\in[\kappa]^\kappa$ for $\N_\alpha$ such that $I_\alpha\subset I_\beta$ for every $\beta<\alpha$. Now by section 2.3 in \cite{SharpeWelch} we get a structure $\P_\alpha$ with $\N_\alpha\in\P_\alpha$ and a $\P_\alpha$-measure $\tilde\mu_\alpha$ on $\kappa$, generated by $I_\alpha$.\footnote{By \textit{generated} here we mean that $X\in\tilde\mu_\alpha$ iff $X$ contains a tail of indiscernibles from $I_\alpha$.} Set $\mu_\alpha:=\tilde\mu_\alpha\cap\M_\alpha$ and let player II play $\mu_\alpha$ in $\G_\zeta(\kappa)$.

	\qquad As the $\mu_\alpha$'s are generated by the $I_\alpha$'s, the $\mu_\alpha$'s are $\subset$-increasing. We have thus created a strategy for player II in $\G_\zeta(\kappa)$ which does not lose at any round $\alpha<\gamma$, making $\kappa$ coherent ${<}\gamma$-Ramsey.
}

The following result is then a direct corollary of Theorems \ref{theo.stratvery} and \ref{theo.succverystrat}.

\qcoro[N.]{
	\label{coro.verystrat}
	For limit ordinals $\alpha$, $\kappa$ is $\omega\alpha$-very Ramsey iff it is coherent ${<}\omega\alpha$-Ramsey. In particular, $\kappa$ is $\lambda$-very Ramsey iff it is strategic $\lambda$-Ramsey for any $\lambda$ with uncountable cofinality.
}

We can now use this equivalence to transfer results from the $\alpha$-very Ramseys over to the strategic versions. The \textit{completely Ramsey cardinals} are the cardinals topping the hierarchy defined in \cite{Feng}. A completely Ramsey cardinal implies the consistency of a Ramsey cardinal, see e.g. Theorem 3.51 in \cite{SharpeWelch}. We are going to use the following characterisation of the completely Ramsey cardinals, which is Lemma 3.49 in \cite{SharpeWelch}.

\qtheo[Sharpe-W.]{
	\label{theo.comvery}
	A cardinal is completely Ramsey if and only if it is $\omega$-very Ramsey.
}

This, together with Theorem \ref{theo.plusone}, immediately yields the following strengthening of Theorem \ref{theo.ramlimram}.

\qcoro[N.]{
	Every $(\omega{+}1)$-Ramsey cardinal is a completely Ramsey stationary limit of completely Ramsey cardinals.
}

The above Theorem \ref{theo.succverystrat} also yields the following consequence.

\coro[N.]{
	Every completely Ramsey cardinal is completely ineffable.
}
\proof{
	From Theorem \ref{theo.comvery} we have that being completely Ramsey is equivalent to being $\omega$-very Ramsey, so the above Theorem \ref{theo.succverystrat} then yields that a completely Ramsey cardinal is coherent ${<}\omega$-Ramsey, which we saw in Theorem \ref{theo.ineff} is equivalent to being completely ineffable.
}

Now, moving to the uncountable case, Corollary \ref{coro.verystrat} yields that strategic $\omega_1$-Ramsey cardinals are $\omega_1$-very Ramsey, and Theorem 3.50 in \cite{SharpeWelch} states that $\omega_1$-very Ramseys are measurable in the core model $K$, assuming $0^\pistol$ doesn't exist, which then shows the following theorem. We also include the original direct proof of that theorem, due to Welch.

\theo[W.]{
	\label{theo.stratramsey}
	Assuming $0^\pistol$ doesn't exist, every strategic $\omega_1$-Ramsey cardinal is measurable in $K$.
}
\proof{
	Let $\kappa$ be strategic $\omega_1$-Ramsey, say $\tau$ is the winning strategy for player II in $\G_{\omega_1}(\kappa)$. Jump to $V[g]$, where $g\subset\col(\omega_1,\kappa^+)$ is $V$-generic. Since $\col(\omega_1,\kappa^+)$ is $\omega$-closed, $V$ and $V[g]$ have the same countable sequences of $V$, so $\tau$ is still a strategy for player II in $\G_{\omega_1}(\kappa)^{V[g]}$, as long as player I only plays elements of $V$.
	
	\qquad Now let $\bra{\kappa_\alpha\mid\alpha<\omega_1}$ be an increasing sequence of regular $K$-cardinals cofinal in $\kappa^+$, let player I in $\G_{\omega_1}(\kappa)$ play $\M_\alpha:=\hull^{H_\theta}(K\l\kappa_\alpha)\prec H_\theta$ and player II follow $\tau$. This results in a countably complete weakly amenable $K$-measure $\mu_{\omega_1}$, which the ``beaver argument''\footnote{See Lemmata 7.3.7--7.3.9 and 8.3.4 in \cite{Zeman} for this argument.} then shows is actually an element of $K$, making $\kappa$ measurable in $K$.
}

A natural question is whether this behaviour persists when going to larger core models. It turns out that the answer is affirmative: every strategic $\omega_1$-Ramsey cardinal is also measurable in Steel's core model below a Woodin, a result due to Schindler which we include with his permission here. We will need the following special case of Corollary 3.1 from \cite{SchindlerIterates}.\footnote{That paper assumes the existence of a measurable as well, but by \cite{JensenSteel} we can omit that here.}

\qtheo[Schindler]{
  \label{theo.Schindler}
  Assume that there exists no inner model with a Woodin cardinal, let $\mu$ be a measure on a cardinal $\kappa$, and let $\pi:V\to\ult(V,\mu)\cong N$ be the ultrapower embedding. Assume that $N$ is closed under countable sequences. Write $K^N$ for the core model constructed inside $N$. Then $K^N$ is a normal iterate of $K$, i.e. there is a normal iteration tree $\T$ on $K$ of successor length such that $\M_{\infty}^{\T}=K^N$. Moreover, we have that $\pi_{0\infty}^{\T}=\pi\restr K$.
}

\theo[Schindler]{
	\label{theo.Schindlerstrat}
  Assuming there exists no inner model with a Woodin cardinal, every strategic $\omega_1$-Ramsey cardinal is measurable in $K$.
}
\proof{
	Fix a large regular $\theta\gg 2^\kappa$. Let $\kappa$ be strategic $\omega_1$-Ramsey and fix a winning strategy $\sigma$ for player II in $\G_{\omega_1}(\kappa)$. Let $g\subset\col(\omega_1,2^\kappa)$ be $V$-generic and in $V[g]$ fix an elementary chain $\bra{M_\alpha\mid\alpha<\omega_1}$ of weak $\kappa$-models $M_\alpha\prec H_\theta^V$ such that $M_\alpha\in V$, $^\omega M_\alpha\subset M_{\alpha+1}$ and $H_{\kappa^+}^V\subset M_{\omega_1}:=\bigcup_{\alpha<\omega_1}M_\alpha$.

	\qquad Note that $V$ and $V[g]$ have the same countable sequences since $\col(\omega_1,2^\kappa)$ is ${<}\omega_1$-closed, so we can apply $\sigma$ to the $M_\alpha$'s, resulting in an $M_{\omega_1}$-measure $\mu$ on $\kappa$. Let $j:M_{\omega_1}\to\ult(M_{\omega_1},\mu)$ be the ultrapower embedding. Since we required that $^\omega M_\alpha\subset M_{\alpha+1}$ we get that $\M_{\omega_1}$ is closed under $\omega$-sequences in $V[g]$, making $\mu$ countably complete in $V[g]$. As we also ensured that $H_{\kappa^+}^V\subset\M_{\omega_1}$ we can lift $j$ to an ultrapower embedding $\pi:V\to\ult(V,\mu)\cong N$ with $N$ transitive.
	
	\qquad Since $V$ is closed under $\omega$-sequences in $V[g]$ we get by standard arguments that $N$ is as well, which means that Theorem \ref{theo.Schindler} applies, meaning that $\pi\restr K:K\to K^N$ is an iteration map with critical point $\kappa$, making $\kappa$ measurable in $K$.
}

\section{Questions and answers}

In this section we give an update on previously posed open questions in the area, as well as posing further open questions. We provide answers for the following questions, which were posed in \cite{HolySchlicht}.

\begin{enumerate}
	\item If $\gamma$ is an uncountable cardinal and the challenger does not have a winning strategy in the game $\G_\gamma^\theta(\kappa)$, does it follow that the judge has one?
	\item If $\omega\leq\alpha\leq\kappa$, are $\alpha$-Ramsey cardinals downwards absolute to the Dodd-Jensen core model?
	\item Does $2$-iterability imply $\omega$-Ramseyness, or conversely?
	\item Does $\kappa$ having the strategic $\kappa$-filter property have the consistency strength of a measurable cardinal?\\
\end{enumerate}

Here the ``challenger'' is player I and the ``judge'' is player II, so this is asking if every $\gamma$-Ramsey is strategic $\gamma$-Ramsey, when $\gamma$ is an uncountable cardinal. Theorem \ref{theo.stratramsey} therefore gives a negative answer to (i) for all uncountable ordinals $\gamma$. Theorem \ref{theo.downK} and Corollary \ref{coro.downK} answer (ii) positively, for $\alpha$-Ramseys with $\alpha$ having uncountable cofinality, and for ${<}\alpha$-Ramseys when $\alpha$ is a limit of limit ordinals. Note that (ii) in the $\alpha=\omega$ case was answered positively in \cite{HolySchlicht}.

\qquad As for (iii), it's mentioned in \cite{HolySchlicht} that Gitman has showed that $\omega$-Ramseys are not in general $2$-iterable by showing that $2$-iterables have strictly stronger consistency strength than the $\omega$-Ramseys, which also follows from Theorem \ref{theo.remlimram} and Theorem 4.8 in \cite{Ramsey2}. Corollary \ref{coro.ind} shows that $\omega$-Ramsey cardinals are $\Delta^2_0$-indescribable, and as $2$-iterables are (at least) $\Pi^1_3$-definable it holds that any $2$-iterable $\omega$-Ramsey cardinal is a limit of $2$-iterables, so that in general $2$-iterables can't be $\omega$-Ramsey either, answering (iii) in the negative. Lastly, Theorem \ref{theo.omegaplusone} gives a positive answer to (iv).

\ques{
	It's not too hard to see that, for a regular uncountable $\lambda$, $\kappa$ is strategic $\lambda$-Ramsey iff there's a ${<}\lambda$-closed forcing $\mathbb P$ such that, in $V^{\mathbb P}$, there's a weakly amenable measure on $\kappa$ with a wellfounded ultrapower. Can we get similar characterisations of strategic $\alpha$-Ramseys for $\alpha$ countable? The proofs of Theorems \ref{theo.virtstrat} and \ref{theo.omegaplusone} give plausible candidates.
}

\ques{
	Are genuine $n$-Ramsey cardinals limits of $n$-Ramsey cardinals? We conjecture this to be true, in analogy with the weakly ineffables being limits of weakly compacts. Since ``weakly ineffable = $\Pi^1_1$-indescribability + subtlety'', this might involve some notion of ``$n$-iterated subtlety''. The difference here is that $n$-Ramseys cannot be \textit{equivalent} to $\Pi^1_{2n+1}$-indescribables for consistency reasons, so there is some work to be done.
}

\ques{
	\label{ques.ctbltheta}
	Fix some $\gamma$ with countable cofinality and an uncountable $\kappa=\kappa^{<\kappa}$. For $\theta>\kappa$ say that $\kappa$ is \textbf{$(\gamma,\theta)$-Ramsey} if player I has no winning strategy in $\G_{\gamma}^\theta(\kappa)$, so that $\kappa$ is $\gamma$-Ramsey iff it's $(\gamma,\theta)$-Ramsey for every $\theta>\kappa$. Do the $(\gamma,\theta)$-Ramseys then eventually form a strict hierarchy? I.e. is there some $\theta>\kappa$ such that $\zfc+\godel{\text{there exists a $(\gamma,\theta_1)$-Ramsey cardinal}}\proves\godel{\text{there exists a $(\gamma,\theta_0)$-Ramsey cardinal}}$ holds for every $\theta_1>\theta_0\geq\theta$? Or, at the opposite end of the spectrum, do the $(\gamma,\theta)$-Ramseys become eventually equivalent? I.e. is there a $\theta>\kappa$ such that $\kappa$ is $(\gamma,\theta_0)$-Ramsey iff it's $(\gamma,\theta_1)$-Ramsey, for all $\theta_1,\theta_0\geq\theta$?
}

\subsection*{Acknowledgements}

The first author would like to thank Ralf Schindler, Philipp Schlicht and Peter Holy for their insightful comments, suggestions and corrections, and Yair Hayut for patiently explaining the above-mentioned Silver-Solovay result on the cut and choose game, which ultimately led to the proof of Theorem \ref{theo.omegaplusone}. We also thank the anonymous referee for their many helpful comments and corrections.

\section{Diagrams}

\subsection*{Consistency implications\footnote{Here dashed lines represent consistency implications which might be equiconsistencies.}}
\begin{center}
	\hspace*{-0.3cm}\includegraphics[scale=0.40]{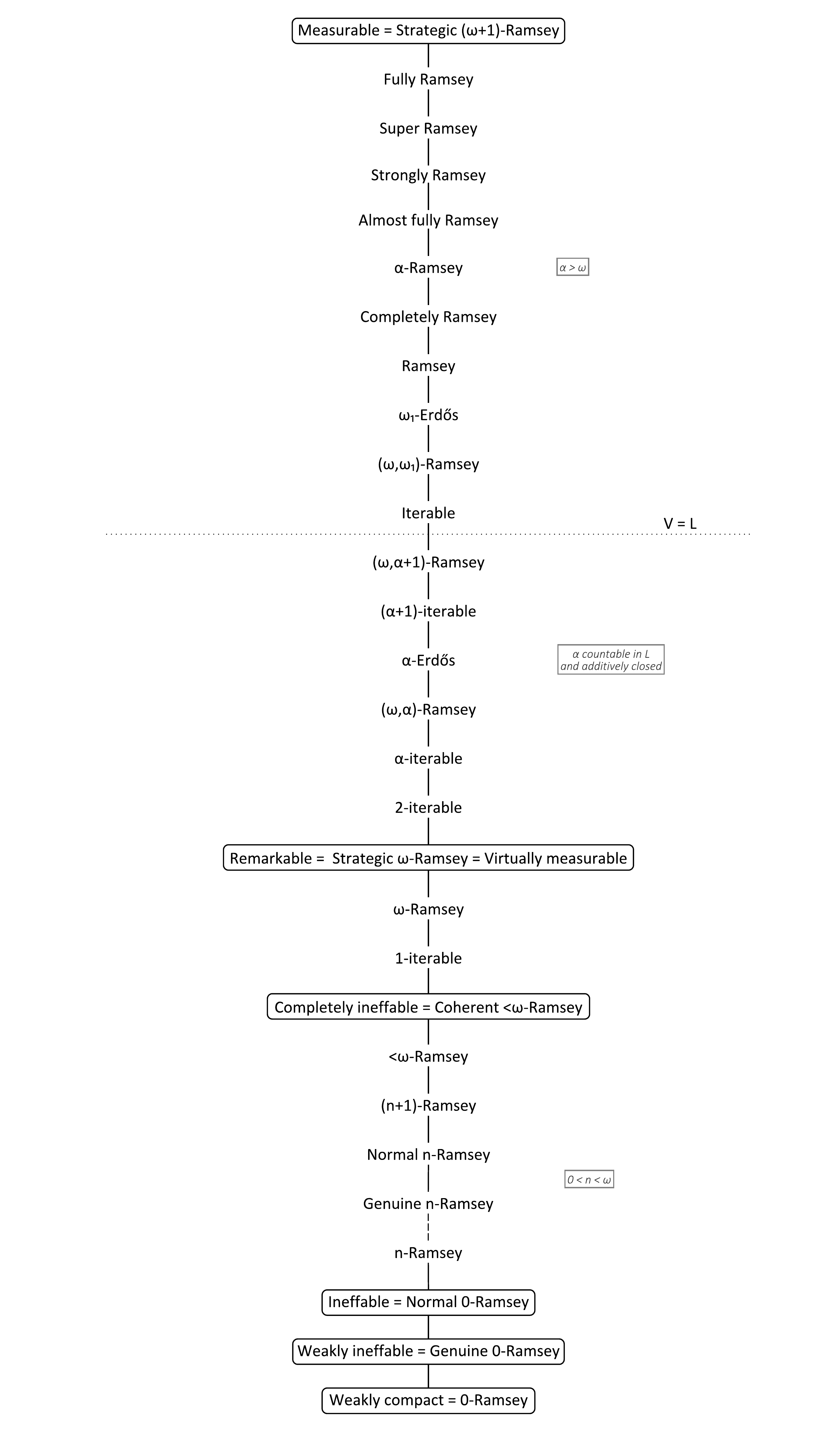}
\end{center}

\pagebreak
\subsection*{Direct implications\footnote{Here dashed lines represent provable direct implications which might be equivalences.}}
\begin{center}
	\hspace*{-2cm}\includegraphics[scale=0.55]{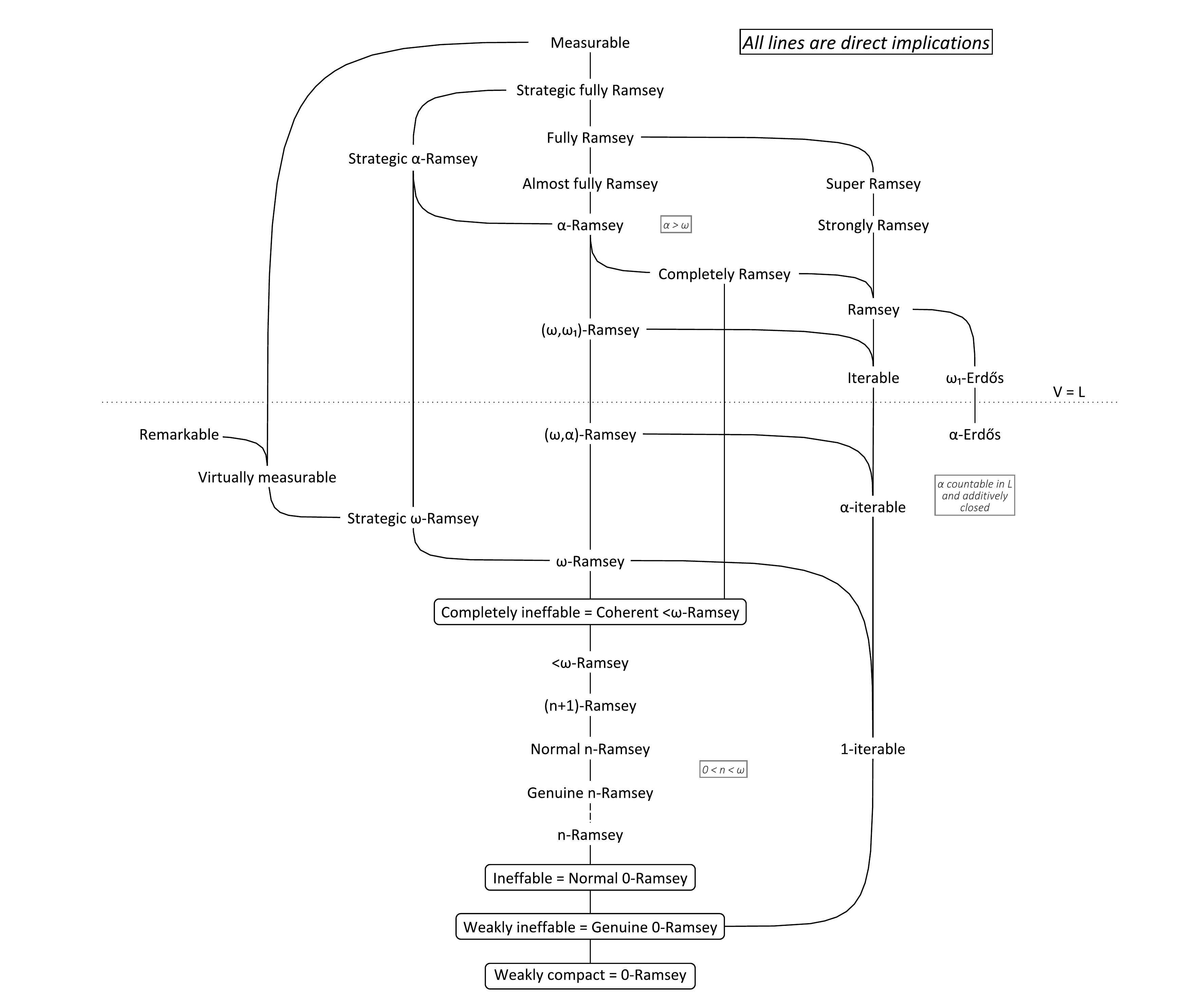}
\end{center}

\pagebreak
\bibliographystyle{apalike}
\bibliography{bib}

\hspace{1em}

\textsc{University of Bristol}

\textsc{Tyndall Ave}

\textsc{Bristol BS8 1TH, UK}

\textit{E-mails:} dan.nielsen@bristol.ac.uk \& p.welch@bristol.ac.uk

\end{document}